\documentclass[dspr]{apjrnl}
\usepackage{amssymb}
\usepackage{amsmath}
\usepackage{amsfonts}
\usepackage[doublespacing]{setspace}

{\theoremstyle{plain}

{\theoremstyle{remark}

}
{\theoremstyle{definition}

}

\iftwocol\AtEndDocument{\end{multicols}}\fi
\input{tcilatex}
\begin{document}

\title{Uniqueness in an Integral Geometry Problem and an Inverse Problem for
the Kinetic Equation}
\author{Arif Amirov$^{\ref{first}}$, Fikret G\"{o}lgeleyen$^{\ref{first}}$
and Masahiro Yamamoto$^{\ref{second}}$ \\
%EndAName
\additem[first]{Department of Mathematics, Bülent Ecevit University,
Zonguldak 67100 Turkey} 
\additem[second]{Department of Mathematical Sciences, 
The University of Tokyo, 
 3-8-1 Komaba, Meguro, Tokyo 153-8914 Japan } \\
E-mail: f.golgeleyen@beun.edu.tr, myama@ms.u-tokyo.ac.jp}
\maketitle

\begin{abstract}
In this paper, we discuss the uniqueness in an integral geometry problem in
a strongly convex domain. Our problem is related to the problem of finding a
Riemannian metric by the distances between all pairs of the boundary points.
For the proof, the problem is reduced to an inverse source problem for a
kinetic equation on a Riemannian manifold and then the uniqueness theorem is
proved in semi-geodesic coordinates by using the tools of Fourier analysis.
\end{abstract}

\section{Introduction}

Let $D\subset 
%TCIMACRO{\U{211d} }%
%BeginExpansion
\mathbb{R}
%EndExpansion
^{n},$ $n\geq 2$ be a simply connected, closed and bounded domain with
boundary $\partial D$ of class $C^{5}$. We assume that the domain $D$ is
strongly convex with respect to a metric $g\in C^{6}(D)$, which means that
for any $x,y\in D$ there exists a unique geodesic $\Gamma (x,y)$ of metric $%
g $ which connects $x$, $y\in D$ and lies in $D$ (e.g., \cite{Cha}).

Henceforth we use the following notations:

$\overset{.}{\xi }=\frac{d}{dt}\xi ,$ $\xi ^{\prime }=(\xi ^{2},...,\xi
^{n}) $ for $\xi =(\xi ^{1},...,\xi ^{n})\in 
%TCIMACRO{\U{211d} }%
%BeginExpansion
\mathbb{R}
%EndExpansion
^{n}$ and $\Gamma (x,y)=\xi (x,y,t)=\left\{ \xi ^{1}(x,y,t),...,\xi
^{n}(x,y,t)\right\} $ is a coordinate representation of the geodesic $\Gamma
(x,y)$. Here $t$ is the natural parameter (see \cite{Dub}, \cite{Rasc}),
that is, $t=Al+B,$ where $l$ is the length of the geodesic (in the metric $g$%
) traced out from some point, $A$ and $B$ are some constants.

We consider the following integral geometry problem:

Throughout this paper, we assume that $a_{ij}=a_{ji}$, $2\leq i,j\leq n$.

\textbf{Problem 1} \textit{Determine the functions} $a_{ij}$\textit{\ in }$D$%
\textit{\ }$(2\leq i,j\leq n)$ \textit{from the integrals}%
\begin{equation}
\int_{\Gamma (x,y)}\left( \underset{i,j=2}{\overset{n}{\sum }}a_{ij}(\xi
(x,y,t))\overset{.}{\xi }^{i}(x,y,t)\overset{.}{\xi }^{j}(x,y,t)\right) dt, 
\tag{1.1}
\end{equation}%
\textit{which are known} \textit{for each pair of the points }$(x,y)\in
\partial D\times \partial D$.

In this paper, we investigate the uniqueness of solution of Problem 1. Our
method for proving the main result, which is stated in Section 2, relies on
the reduction of Problem 1 to some kinetic equation (see (3.5) below) on a
Riemannian manifold where the metric is considered in a semi-geodesic system
of coordinates.

Here we note that when we replace $i,j=2$ in the summation by $i,j=1$ in
(1.1), we do not know the uniqueness in determining $a_{ij}$, $1\leq i,j\leq
n.$ The choice of the indices "$i,j\neq 1$" depends on the semi-geodesic
coordinates which we use in this paper. For detailed explanations, see also
Remark 1 in Section 4.

Problem 1 is related to an inverse problem of determining the Riemannian
metric by the distances between boundary points, see Chapter 1 of \cite{Sha}%
; also \cite{A-dok}, \cite{A-tok}. Such an inverse problem is the
mathematical model of several important medical imaging techniques and
geophysical problems, and has called wide attention. As for the uniqueness
theorem and stability estimates, see Muhometov \cite{Mukh} in two
dimensions, and a recent work Pestov and Uhlmann \cite{PesUh}. For higher
dimensions, we refer to Bernstein and Gerver \cite{Berst}, Beylkin \cite%
{Beyl}, Muhometov and Romanov \cite{MukhRo}. Also see e.g., \cite{Natt}, 
\cite{Patern1}. As for other kinds of inverse problems from the integral
geometry, we refer to \cite{AnRo}, \cite{Dai}, \cite{Gelf}, \cite{Patern2}, 
\cite{Pestov}, \cite{Sha} - \cite{SteU2}. Here we do not intend any complete
lists of references. The connections between these problems and the inverse
problems for parabolic, hyperbolic and kinetic equations are described in
the works \cite{A-bo}, \cite{Ani}, \cite{Klib}, \cite{klib2}, \cite{Lav}, 
\cite{Rom}.

\section{Main Result}

Throughout this paper, we set%
\begin{equation*}
%TCIMACRO{\U{211d} }%
%BeginExpansion
\mathbb{R}
%EndExpansion
_{0}^{n}:=%
%TCIMACRO{\U{211d} }%
%BeginExpansion
\mathbb{R}
%EndExpansion
^{n}\backslash \left\{ 0\right\} ;\text{ }%
%TCIMACRO{\U{211d} }%
%BeginExpansion
\mathbb{R}
%EndExpansion
_{0}^{n-1}:=%
%TCIMACRO{\U{211d} }%
%BeginExpansion
\mathbb{R}
%EndExpansion
^{n-1}\backslash \left\{ 0\right\} \text{,}
\end{equation*}%
\begin{equation*}
C_{0}^{5}(D)=\{a\in C^{5}(D)|\,\mbox{supp}\,a\subset D\}\text{.}
\end{equation*}%
We note that $D$ is closed and supp $a$ is not necessarily a proper subset
of $D$.

The result which we have obtained for Problem 1 is given by Theorem 1:

\textbf{Theorem 1 }\textit{Let }$a_{ij}=a_{ji}$, $2\leq i,j\leq n$, and 
\textit{the domain }$D$\textit{\ be strongly convex with respect to the
metric }$g\in C^{6}(D)$.\textit{\ Then Problem 1 may have at most one
solution }$a_{ij}$ such that $a_{ij}\in C_{0}^{5}(D)$, where $2\leq i,j\leq
n $.

It is worth noting here that, in the proof, we assume that the metric has
the property $g_{11}=1$, $g_{1i}=0$, $2\leq i\leq n,$ which is related to
the semi-geodesic coordinates and extensively used in the theory of
relativity (see, e.g., \cite{Petrov} and the references therein). We more
precisely explain as follows: It is known that if there exists a point $%
x^{0}\in D$ such that any point $x\in D$\ can be joined with $x^{0}$ by a
unique geodesic of the metric $g,$ then the metric $g$ has a semi-geodesic
coordinates in $D$ (see \cite{Rasc}, p. 448). Hence, for any domain $D$
which is strongly convex with respect to the metric $g,$ one can at least
locally introduce a semi-geodesic system of coordinates at any $x\in D$ with
respect to $g$. Moreover, in the semi-geodesic system of coordinates $x^{i},$
we have $g_{11}=1,$ $g_{1i}=0,$ $2\leq i\leq n$, locally. Conversely, these
conditions are sufficient for the system with coordinates $x^{i}$ to be
semi-geodesic for the metric $g$ in $D$ (see \cite{Rasc}, p. 450 and \cite%
{wae}, p. 76).

This paper consists of six sections and one appendix. The rest part of the
paper is organized as follows. In Section 3, we reduce Problem 1 to an
inverse source problem for a kinetic equation on a Riemannian manifold. In
Section 4, we present three lemmata and apply the generalized Fourier
transform to the kinetic equation. In Section 5, we reformulate our problem
by introducing Riemannian coordinates and prove some auxiliary results.
Finally, Section 6 is devoted to the proof of the main result: Theorem 1.
The proofs of the lemmata which are presented in Section 4 are given in
Appendix.

\section{Reduction of Problem 1 to an Inverse Problem}

Let us introduce the function%
\begin{equation}
u(x,\xi )=\overset{n}{\underset{i,j=2}{\sum }}\int_{\gamma (x,\xi
)}a_{ij}\left( z\left( x,\xi ,t\right) \right) \overset{.}{z}^{i}\left(
x,\xi ,t\right) \overset{.}{z}^{j}\left( x,\xi ,t\right) dt,  \tag{3.1}
\end{equation}%
where $\gamma \left( x,\xi \right) $ is the ray of the metric $g=\left(
g_{ij}\right) $ starting from$\ x\in D$ in direction $\xi \in \mathbb{R}_0^n$
and $a_{ij}\in C^{5}_0(D)$.

It is known (see \cite{Carmo}, \cite{Dub}, \cite{Rasc}) that $\gamma \left(
x,\xi \right) =\left( z^{1}\left( x,\xi ,t\right) ,...,z^{n}\left( x,\xi
,t\right) \right) $ is the solution of the following system of differential
equations%
\begin{equation}
\frac{d^{2}z^{i}}{dt^{2}}=-\Gamma _{jk}^{i}(z)\overset{.}{z}^{j}\overset{.}{z%
}^{k},\text{ }1\leq i\leq n\text{,}  \tag{3.2}
\end{equation}%
with the Cauchy data%
\begin{equation}
z(0)=x,\ \overset{.}{z}(0)=\xi ,  \tag{3.3}
\end{equation}%
where $\Gamma _{jk}^{i}$ are the Christoffel symbols of the metric $g$ and%
\begin{eqnarray*}
z\left( x,\xi ,t\right) &=&\left( z^{1}\left( x,\xi ,t\right)
,...,z^{n}\left( x,\xi ,t\right) \right) , \\
\overset{.}{z}\left( x,\xi ,t\right) &=&\left( \overset{.}{z}^{1}\left(
x,\xi ,t\right) ,...,\overset{.}{z}^{n}\left( x,\xi ,t\right) \right) ,\text{
}\overset{.}{z}^{i}=\frac{d}{dt}z^{i}\text{.}
\end{eqnarray*}%
We can prove that the solution of Problem (3.2)-(3.3) has the following
property:%
\begin{equation}
z\left( x,\xi ,t\right) =z\left( x,\nu ,\left\vert \xi \right\vert t\right) ;%
\text{ }\overset{.}{z}\left( x,\xi ,t\right) =\left\vert \xi \right\vert 
\overset{.}{z}\left( x,\nu ,\left\vert \xi \right\vert t\right) ,  \tag{3.4}
\end{equation}%
where $\nu =\xi /\left\vert \xi \right\vert $ and $\left\vert \xi
\right\vert ^{2}=\underset{i,j=1}{\overset{n}{\sum }}g_{ij}\xi ^{i}\xi ^{j}$
(e.g., Lemma 2.6 on p.64 in do Carmo \cite{Carmo}).

Let $G^{\prime }$ denote a closed, bounded set of variables $\xi ^{\prime
}=(\xi ^{2},...,\xi ^{n})$ such that $0\notin G^{\prime }$ and let $%
G=\left\{ \xi \in 
%TCIMACRO{\U{211d} }%
%BeginExpansion
\mathbb{R}
%EndExpansion
^{n}\mbox{ }|\mbox{ }\xi =(\xi ^{1},\xi ^{\prime }),\mbox{ }\xi ^{1}\in 
%TCIMACRO{\U{211d} }%
%BeginExpansion
\mathbb{R}
%EndExpansion
^{1},\ \xi ^{\prime }\in G^{\prime }\right\} ,$ $\Omega =\left\{ (x,\xi )%
%\text{ \TEXTsymbol{\vert} }
|x\in D,\mbox{ }\xi \in G\right\} .$

Differentiating both sides of (3.1) at the point $x$ in the direction $\xi $
and by using (3.2) and (3.3), we have the following kinetic equation%
\begin{equation}
\underset{j=1}{\overset{n}{\mathop{\textstyle \sum }}}\xi ^{j}\frac{\partial
u}{\partial x^{j}}\underset{j,k,s=1}{-\overset{n}{\mathop{\textstyle \sum }}}%
\Gamma _{jk}^{^{s}}(x)\xi ^{k}\xi ^{j}\frac{\partial u}{\partial \xi ^{s}}=%
\underset{j,k=2}{\overset{n}{\mathop{\textstyle \sum }}}a_{jk}(x)\xi ^{k}\xi
^{j}.  \tag{3.5}
\end{equation}%
By the setting of Problem 1, using equalities (3.1), (3.4) and the fact that 
$a_{jk}(x)$ is zero outside of $D,$ we conclude that the function $u(x,\xi )$
is known for $(x,\xi )\in \partial D\times 
%TCIMACRO{\U{211d} }%
%BeginExpansion
\mathbb{R}
%EndExpansion
_{0}^{n}$. Then it is easy to see that the uniqueness of the solution of
Problem 1 follows from the uniqueness of the solution of the following
problem:

\textbf{Problem 2 }\textit{Determine a matrix-valued function }$\left(
a_{jk}\right) $, $(2\leq j,k\leq n)$ \textit{from equation (3.5) provided
that }$u(x,\xi )$\textit{\ is known for} $(x,\xi )\in \partial D\times 
%TCIMACRO{\U{211d} }%
%BeginExpansion
\mathbb{R}
%EndExpansion
_{0}^{n}$\textit{.}

In order to prove the uniqueness for Problem 2, it is sufficient to assume%
\begin{equation}
u(x,\xi )=0\text{,}\ (x,\xi )\in \partial D\times 
%TCIMACRO{\U{211d} }%
%BeginExpansion
\mathbb{R}
%EndExpansion
_{0}^{n}.  \tag{3.6}
\end{equation}

\section{Application of the Fourier Transform to Problem 2}

In this section, we present three lemmata which describe some important
properties of the function $u(x,\xi )$ and its Fourier transform. The proofs
are given in Appendix. The Fourier transform $\mathcal{F}(I)$ of a function $%
I(x,\xi )\in L_{1}\left( 
%TCIMACRO{\U{211d} }%
%BeginExpansion
\mathbb{R}
%EndExpansion
_{\xi ^{1}}^{1}\right) $\ with respect to the variable$\ \xi ^{1}$ is
defined by%
\begin{equation*}
\mathcal{F}\left\{ I(x,\xi )\right\} =\widehat{I}\left( x,\eta ,\xi ^{\prime
}\right) :=\int_{-\infty }^{\infty }I\left( x,\xi ^{1},\xi ^{\prime }\right)
e^{-\sqrt{-1}\xi ^{1}\eta }d\xi ^{1}\text{,}
\end{equation*}%
where $\eta $\ is the dual variable of $\xi ^{1}$. Henceforth we set%
\begin{eqnarray*}
\partial _{x^{s}} &=&\dfrac{\partial }{\partial x_{s}},\text{ }\partial
_{\eta }=\dfrac{\partial }{\partial \eta },\text{ }\partial _{\xi ^{j}}=%
\dfrac{\partial }{\partial \xi ^{j}}, \\
\partial _{\xi }^{\beta } &=&\partial _{\xi ^{1}}^{\beta _{1}}\cdots
\partial _{\xi ^{n}}^{\beta _{n}},\mathit{\ }\left\vert \beta \right\vert
=\beta _{1}+...+\beta _{n}, \\
\partial _{\xi ^{\prime }}^{\beta ^{\prime }} &=&\partial _{\xi ^{2}}^{\beta
_{2}}\cdots \partial _{\xi ^{n}}^{\beta _{n}},\mathit{\ }\left\vert \beta
^{\prime }\right\vert =\beta _{2}+...+\beta _{n}.
\end{eqnarray*}%
Moreover we introduce the sets%
\begin{eqnarray*}
\Delta _{\eta }^{\rho } &:&=\left\{ \eta \in 
%TCIMACRO{\U{211d} }%
%BeginExpansion
\mathbb{R}
%EndExpansion
_{\eta }^{1}|\mbox{ }\rho \eta >0\right\} , \\
\overline{\Delta }_{\eta }^{\rho } &:&=\left\{ \eta \in 
%TCIMACRO{\U{211d} }%
%BeginExpansion
\mathbb{R}
%EndExpansion
_{\eta }^{1}|\mbox{ }\rho \eta \geq 0\right\} ,
\end{eqnarray*}%
where we set $\rho =-1,1$.

Let us introduce the auxiliary functions%
\begin{equation}
I_{ij}(x,\xi )=\int_{\gamma (x,\xi )}b\left( z\left( x,\xi ,t\right) \right) 
\overset{.}{z}^{i}\left( x,\xi ,t\right) \overset{.}{z}^{j}\left( x,\xi
,t\right) dt,\text{ }2\leq i,j\leq n,  \tag{4.1}
\end{equation}%
where the function $b\in C^{5}(%
%TCIMACRO{\U{211d} }%
%BeginExpansion
\mathbb{R}
%EndExpansion
^{n})$ is zero outside of $D$. Note that the function $u(x,\xi )$ in (3.1)
is defined as the sum of the functions of the form $I_{ij}$. By taking into
account (3.4), let us rewrite (4.1) in the form%
\begin{eqnarray}
I_{ij}(x,\xi ) &=&\int_{0}^{\infty }b\left( z\left( x,\xi ,t\right) \right) 
\overset{.}{z}^{i}\left( x,\xi ,t\right) \overset{.}{z}^{j}\left( x,\xi
,t\right) dt  \notag \\
&=&\int_{0}^{\infty }b\left( z\left( x,\nu ,\left\vert \xi \right\vert
t\right) \right) \left\vert \xi \right\vert \overset{.}{z}^{i}\left( x,\nu
,\left\vert \xi \right\vert t\right) \left\vert \xi \right\vert \overset{.}{z%
}^{j}\left( x,\nu ,\left\vert \xi \right\vert t\right) dt  \notag \\
&=&\frac{1}{\left\vert \xi \right\vert }\int_{0}^{\infty }b\left( z\left(
x,\nu ,\tau \right) \right) \left\vert \xi \right\vert \overset{.}{z}%
^{i}\left( x,\nu ,\tau \right) \left\vert \xi \right\vert \overset{.}{z}%
^{j}\left( x,\nu ,\tau \right) d\tau \text{.}  \TCItag{4.2}
\end{eqnarray}%
\textbf{Lemma 1 }\textit{Let }$D$\textit{\ be strongly convex with respect
to the metric }$g=(g_{ij})\in C^{6}\left( D\right) $.\textit{\ Then the
functions }$I_{ij}$ satisfy the following properties:\textit{\ }

(i) $\partial _{\xi }^{\beta }I_{ij},$\textit{\ }$\partial _{\xi }^{\beta
}\partial _{x^{s}}I_{ij}\in C\left( \Omega \right) $\textit{\ for }$0\leq
\left\vert \beta \right\vert \leq 4,$

(ii) \textit{For fixed }$x\in D$\textit{\ and }$\xi ^{\prime }\in G^{\prime
} $\textit{\ }$\left( \xi ^{\prime }\neq 0\right) ,$

\quad (a) $\partial _{\xi }^{\beta }I_{ij}$\textit{,}$\ \partial _{\xi
}^{\beta }\partial _{x^{s}}I_{ij}\in L_{2}\left( 
%TCIMACRO{\U{211d} }%
%BeginExpansion
\mathbb{R}
%EndExpansion
_{\xi ^{1}}^{1}\right) $ \textit{for\ }$\left\vert \beta \right\vert \leq 2,$

\quad (b) $\partial _{\xi }^{\beta }I_{ij}$\textit{,}$\ \partial _{\xi
}^{\beta }\partial _{x^{s}}I_{ij}\in L_{1}\left( 
%TCIMACRO{\U{211d} }%
%BeginExpansion
\mathbb{R}
%EndExpansion
_{\xi ^{1}}^{1}\right) \cap L_{2}\left( 
%TCIMACRO{\U{211d} }%
%BeginExpansion
\mathbb{R}
%EndExpansion
_{\xi ^{1}}^{1}\right) $\textit{\ for }$\left\vert \beta \right\vert =3,$

\quad (c) $\xi ^{1}\partial _{\xi }^{\beta }I_{ij}$\textit{,}$\ \xi
^{1}\partial _{\xi }^{\beta }\partial _{x^{s}}I_{ij}\in L_{1}\left( 
%TCIMACRO{\U{211d} }%
%BeginExpansion
\mathbb{R}
%EndExpansion
_{\xi ^{1}}^{1}\right) \cap L_{2}\left( 
%TCIMACRO{\U{211d} }%
%BeginExpansion
\mathbb{R}
%EndExpansion
_{\xi ^{1}}^{1}\right) $\textit{\ for }$\left\vert \beta \right\vert =4,$%
\textit{\newline
where }$1\leq s\leq n$ and $2 \le i,j \le n$.

\textbf{Remark 1. }It is worth to note that Lemma 1 is not valid if at least
one of the indices $i$ or $j$ in (4.2)\ is equal to $1$.

Indeed we note:

(i) Since $\overset{.}{z}^{1}(x,\nu ,t)\rightarrow \overset{.}{z}^{1}(x,\pm
\nu ^{0},t)=\pm 1$ as $\xi ^{1}\rightarrow \pm \infty ,$ the function $%
\left\vert \xi ^{1}\overset{.}{z}^{1}(x,\nu ,t)\right\vert $ increases like $%
\left\vert \xi ^{1}\right\vert $ as $\xi ^{1}\rightarrow \infty $,

(ii) Integral (4.2) is taken on a finite interval $\left[ 0,d_{0}\right] $,
where $d_{0}$ is the diameter of $D$ in the metric $g=\left(
g_{ij}(x)\right) $,

(iii) $\left\vert \mbox{ }\left\vert \xi \right\vert \mbox{ }\overset{.}{z}%
^{k}\left( x,\nu ,t\right) \right\vert \leq K_{1},$ (See the proof of Lemma
1 in Appendix).

Then it follows from (4.2) that if only one of the indices $i,$ $j$ is equal
to 1, then$\ I_{ij}(x,\xi )$ is bounded only on the set $\Omega ,$ whereas
the function $\left\vert I_{11} (x,\xi )\right\vert $ increases like $%
\left\vert \xi^{1}\right\vert $ as $\xi ^{1}\rightarrow \infty$ for each
fixed $(x,\xi ^{\prime })\in D\times G^{\prime }$.

\textbf{Lemma 2 }\textit{Let the conditions of Lemma 1 be satisfied. Then we
have}

(i) $\partial _{\xi ^{\prime }}^{\beta ^{\prime }}\widehat{I}_{ij},$\textit{%
\ }$\partial _{\xi ^{\prime }}^{\beta ^{\prime }}\partial _{x^{_{s}}}%
\widehat{I}_{ij}\in C\left( D\times \Delta _{\eta }^{\rho }\times G^{\prime
}\right) \cap L_{2}\left( 
%TCIMACRO{\U{211d} }%
%BeginExpansion
\mathbb{R}
%EndExpansion
_{\eta }^{1}\right) $\textit{\ for }$\left\vert \beta ^{\prime }\right\vert
\leq 2,$

(ii) $\partial _{\xi ^{\prime }}^{\beta ^{\prime }}\widehat{I}_{ij},$\textit{%
\ }$\partial _{\xi ^{\prime }}^{\beta ^{\prime }}\partial _{x^{_{s}}}%
\widehat{I}_{ij}\in C\left( D\times 
%TCIMACRO{\U{211d} }%
%BeginExpansion
\mathbb{R}
%EndExpansion
_{\eta }^{1}\times G^{\prime }\right) \cap L_{2}\left( 
%TCIMACRO{\U{211d} }%
%BeginExpansion
\mathbb{R}
%EndExpansion
_{\eta }^{1}\right) $\textit{\ for }$\left\vert \beta ^{\prime }\right\vert
=3,$

(iii) $\eta ^{r}\partial _{\xi ^{\prime }}^{\beta ^{\prime }}\partial _{\eta
}\widehat{I}_{ij},\ \eta ^{r}\partial _{\xi ^{\prime }}^{\beta ^{\prime
}}\partial _{\eta }\partial _{x^{_{s}}}\widehat{I}_{ij}\in C\left( D\times 
%TCIMACRO{\U{211d} }%
%BeginExpansion
\mathbb{R}
%EndExpansion
_{\eta }^{1}\times G^{\prime }\right) \cap L_{2}\left( 
%TCIMACRO{\U{211d} }%
%BeginExpansion
\mathbb{R}
%EndExpansion
_{\eta }^{1}\right) $\textit{\ for }$r+\left\vert \beta ^{\prime
}\right\vert =4,$\textit{\ }$0\leq r\leq 4$,

\textit{where} $1\leq s\leq n$ and $2\leq i,j\leq n$.

\textbf{Remark 2.} Since the function $u(x,\xi )$ is defined as the sum of
the functions of the form $I_{ij}$, Lemmata 1 and 2 are also valid for $%
u(x,\xi )$. On the other hand, due to the reason mentioned in Remark 1, if
the right-hand side of equality (3.1) is replaced by the expression%
\begin{equation*}
\overset{n}{\underset{i,j=1}{\sum }}\int_{\gamma (x,\xi )}a_{ij}\left(
z\left( x,\xi ,t\right) \right) \overset{.}{z}^{i}\left( x,\xi ,t\right) 
\overset{.}{z}^{j}\left( x,\xi ,t\right) dt,
\end{equation*}%
then these lemmata are not valid.

Hence, for each fixed $x\in D$ and $\xi ^{\prime }\in G^{\prime },$ it is
possible to apply the generalized Fourier transform with respect to variable 
$\xi ^{1}$ to equation (3.5). Then we have%
\begin{multline}
\sqrt{-1}\partial _{_{\eta }}\partial _{x^{1}}\widehat{u}-2\sqrt{-1}\underset%
{j,k=2}{\overset{n}{\sum }}\Gamma _{1k}^{j}\xi ^{k}\partial _{_{\eta
}}\partial _{\xi ^{j}}\widehat{u}+\underset{j=2}{\overset{n}{\sum }}\xi
^{j}\partial _{x^{_{j}}}\widehat{u}  \notag \\
-\sqrt{-1}\underset{j,k=2}{\overset{n}{\sum }}\Gamma _{jk}^{1}\xi ^{k}\xi
^{j}\eta \widehat{u}-\underset{j,k,s=2}{\overset{n}{\sum }}\Gamma
_{jk}^{s}\xi ^{k}\xi ^{j}\partial _{\xi ^{s}}\widehat{u}=2\pi \delta \left(
\eta \right) \underset{k,j=2}{\overset{n}{\sum }}a_{kj}\left( x\right) \xi
^{k}\xi ^{j},  \tag{4.3}
\end{multline}%
where we use the fact that $\Gamma _{1k}^{1}=\Gamma _{11}^{k}=0,$ $1\leq
k\leq n$. In (4.3), we recall that $\delta \left( \eta \right) $ is the
Dirac delta function and $\mathcal{F}(1)=2\pi \delta \left( \eta \right) $.

By Remark 2 and taking into account Lemma 2, we see that the functions $%
\widehat{u}$, $\partial _{\eta }\widehat{u}$ are continuously differentiable
in the region $D\times G^{\prime }$ for both cases $\eta >0$ and $\eta <0$.
Then using (4.3), we obtain that for both cases $\eta >0$ and $\eta <0$ the
function$\ \widehat{u}$ satisfies the equation%
\begin{gather}
\sqrt{-1}\partial _{_{\eta }}\partial _{x^{1}}\widehat{u}-2\sqrt{-1}\underset%
{j,k=2}{\overset{n}{\sum }}\Gamma _{1k}^{j}\xi ^{k}\partial _{_{\eta
}}\partial _{\xi ^{j}}\widehat{u}+\underset{j=2}{\overset{n}{\sum }}\xi
^{j}\partial _{x^{_{j}}}\widehat{u}  \notag \\
-\sqrt{-1}\underset{j,k=2}{\overset{n}{\sum }}\Gamma _{jk}^{1}\xi ^{k}\xi
^{j}\eta \widehat{u}-\underset{j,k,s=2}{\overset{n}{\sum }}\Gamma
_{jk}^{s}\xi ^{k}\xi ^{j}\partial _{\xi ^{s}}\widehat{u}=0  \tag{4.4}
\end{gather}%
in the classical sense.

Putting$\ \widehat{u}=p+\sqrt{-1}q$ and separating the real and the
imaginary parts of the left-hand side of equation (4.4) for both$\ \eta >0$
and $\eta <0$, we have the following equations with respect to the functions$%
\ \partial _{\eta }p$ and $\partial _{_{\eta }}q$ respectively:%
\begin{equation}
\partial _{\eta }\partial _{x^{1}}p-2\underset{j,k=2}{\overset{n}{%
\mathop{\textstyle \sum }}}\Gamma _{1k}^{j}\xi ^{k}\partial _{\eta }\partial
_{\xi ^{j}}p=\digamma _{1},  \tag{4.5}
\end{equation}%
\begin{equation}
\partial _{\eta }\partial _{x^{1}}q-2\underset{j,k=2}{\overset{n}{%
\mathop{\textstyle \sum }}}\Gamma _{1k}^{j}\xi ^{k}\partial _{\eta }\partial
_{\xi ^{j}}q=\digamma _{2},  \tag{4.6}
\end{equation}%
where%
\begin{equation*}
\digamma _{1}=\underset{j,k=2}{\overset{n}{\sum }}\Gamma _{jk}^{1}\xi
^{k}\xi ^{j}\eta p-\underset{j=2}{\overset{n}{\mathop{\textstyle \sum }}}\xi
^{j}\partial _{x^{j}}q+\underset{s,j,k=2}{\overset{n}{\sum }}\Gamma
_{jk}^{s}\xi ^{k}\xi ^{j}\partial _{\xi ^{s}}q,
\end{equation*}%
\begin{equation*}
\digamma _{2}=\underset{j,k=2}{\overset{n}{\sum }}\Gamma _{jk}^{1}\xi
^{k}\xi ^{j}\eta q+\underset{j=2}{\overset{n}{\mathop{\textstyle \sum }}}\xi
^{j}\partial _{x^{j}}p-\underset{s,j,k=2}{\overset{n}{\sum }}\Gamma
_{jk}^{s}\xi ^{k}\xi ^{j}\partial _{_{\xi ^{s}}}p.
\end{equation*}%
In this work, the uniqueness of the solution of the problem is investigated
under the assumption of the existence of the solution$.$ Hence it is assumed
that there exists a solution $u(x,\xi )$ to equation (3.5), i.e., there
exists a solution $\widehat{u}=p+\sqrt{-1}q$ to equation (4.3) which
satisfies the properties indicated in Remark 2 and the condition%
\begin{equation}
u\left( x,\xi \right) =0,\text{ }(x,\xi )\in \partial D\times G,\text{ }%
\left( \widehat{u}(x,\eta ,\xi ^{\prime })=0,\ (x,\xi ^{\prime })\in
\partial D\times G^{\prime }\right) \text{.}  \tag{4.7}
\end{equation}

\textbf{Lemma 3. }\textit{Let the conditions of Lemma 1 and condition (3.6)
be satisfied. Then}%
\begin{equation*}
\partial _{\xi ^{\prime }}^{\beta ^{\prime }}\partial _{\eta }\widehat{u}%
(x,\cdot ,\xi ^{\prime }),\,\partial _{\xi ^{\prime }}^{\beta ^{\prime
}}\partial _{x^{s}}\partial _{\eta }\widehat{u}(x,\cdot ,\xi ^{\prime })\in
L_{1}\left( \Delta _{\eta }^{\rho }\right) \cap L_{2}\left( \Delta _{\eta
}^{\rho }\right)
\end{equation*}%
\textit{for fixed} $(x,\xi ^{\prime })\in D\times G^{\prime }$ and 
\begin{eqnarray*}
\partial _{\xi ^{\prime }}^{\beta ^{\prime }}\partial _{\eta }\widehat{u}%
,\,\partial _{\xi ^{\prime }}^{\beta ^{\prime }}\partial _{x^{s}}\partial
_{\eta }\widehat{u} &\in &C\left( D\times \Delta _{\eta }^{\rho }\times
G^{\prime }\right) , \\
\partial _{\xi ^{\prime }}^{\beta ^{\prime }}\widehat{u},\mathit{\ }\partial
_{\xi ^{\prime }}^{\beta ^{\prime }}\partial _{x^{s}}\widehat{u} &\in
&C\left( D\times \overline{\Delta }_{\eta }^{\rho }\times G^{\prime }\right)
,
\end{eqnarray*}%
\textit{for }$0\leq \left\vert \beta ^{\prime }\right\vert \leq 2$ \textit{%
and} $1\leq s\leq n$\textit{.}

\textbf{Remark 3.} In the proof of Lemma 3, we essentially use condition
(4.7).

\section{Some auxiliary results in Riemannian coordinates}

We set $D_{\varepsilon }=\left\{ x\in 
%TCIMACRO{\U{211d} }%
%BeginExpansion
\mathbb{R}
%EndExpansion
^{n}\mbox{ }/\mbox{ }d(x,D)<\varepsilon \right\} $ where $d(x,D)=\inf_{y\in
D}\left\vert x-y\right\vert .$ For some $\varepsilon >0,$\ let $%
D_{\varepsilon }$\ be strongly convex with respect to the metric $g=\left(
g_{ij}^{{}}\right) \in C^{6}(\overline{D_{\varepsilon }}).$ We assume that
the metric $g=\left( g_{ij}^{{}}(x)\right) $ is Euclidean outside of the
region $D_{\varepsilon }$. Then as in \cite{Rasc} (see p. 506) it is
possible to construct the Riemannian coordinates with center $\widetilde{x}%
_{0}\in D$ by the formula $y^{i}=\xi ^{i}t$ $(1\leq i\leq n)$, where $t$ is
the same parameter as in system (3.2). Let us also note that the coordinate $%
y^{i}$ does not depend on the selection of the natural parameter $t$ on this
geodesic (see \cite{Rasc}, p. 506 and \cite{Eis}, p. 53). It is clear that
we have\ $y^{i}=0$ $\left( 1\leq i\leq n\right) $ at the point $\widetilde{x}%
_{0}$. Since $g=\left( g_{ij}\right) \in C^{6}(\overline{D}_{\varepsilon })$%
, the solution $z(x_{0},\xi ,t)$ of Problem (3.2)-(3.3) is of $C^{5}$-class
with respect to $t$ and $\xi $. Then from the determination of the
Riemannian coordinates given above, it follows that $%
y^{i}=y^{i}(x^{1},...,x^{n})\in C^{5}(\overline{D}_{\varepsilon })$\ ($1\leq
i\leq n$). Moreover the $x^{1}$-axis is the same as the $y^{1}$-axis,
because the straight line passing from $\widetilde{x}_{0}$ and parallel to
the $x^{1}$-axis is a geodesic of the metric $g$ which corresponds to the
vector $(1,0,...,0)\in 
%TCIMACRO{\U{211d} }%
%BeginExpansion
\mathbb{R}
%EndExpansion
^{n}$ on which $x^{1}$ is length of the geodesic. Then in the new
coordinates the differential form $\underset{i,j=2}{\overset{n}{\sum }}%
a_{ij}(x)\overset{.}{\xi }^{i}\overset{.}{\xi }^{j}$ takes the form $%
\underset{k,s=2}{\overset{n}{\sum }}\widetilde{a}_{ks}(y)\overset{.}{\zeta }%
^{k}\overset{.}{\zeta }^{s}$ where $\widetilde{a}_{ks}=a_{ij}\partial
_{y^{k}}x^{i}\partial _{y^{s}}x^{j}$ and $a_{ij}=a_{ij}\left( x\left(
y\right) \right) $.

In this case, we introduce the function%
\begin{equation}
\widetilde{u}\left( y,\zeta \right) =\overset{n}{\underset{i,j=2}{\sum }}%
\int_{\widetilde{\gamma }(y,\zeta )}\widetilde{a}_{ij}\left( \tilde{z}\left(
y,\zeta ,t\right) \right) \overset{.}{\tilde{z}}^{i}\left( y,\zeta ,t\right) 
\overset{.}{\tilde{z}}^{j}\left( y,\zeta ,t\right) dt,  \tag{5.1}
\end{equation}%
analoguos to formula (3.1), where $\widetilde{\gamma }(y,\zeta )=\left\{ 
\tilde{z}^{1}\left( y,\zeta ,t\right) ,...,\tilde{z}^{n}\left( y,\zeta
,t\right) \right\} $ is a geodesic of the metric $g$ in the new coordinates
passing from the point $y\in \widetilde{D}$ in direction $\zeta $. It can be
easily proved that Lemmata 1-3 are also valid for the function $\widetilde{u}%
\left( y,\zeta \right) $ by using the similar arguments which we used in the
proofs of Lemmata 1-3 for the function $u(x,\xi )$. Then, by $\widetilde{%
\Gamma }_{jk}^{^{i}}=0$ at $y=0,$ the function $\widetilde{u}\left( y,\zeta
\right) $ satisfies the following analogue of equation (3.5):%
\begin{equation}
\underset{i=1}{\overset{n}{\mathop{\textstyle \sum }}}\zeta ^{i}\partial
_{y^{i}}\widetilde{u}=\underset{s,k=2}{\overset{n}{\mathop{\textstyle \sum }}%
}\widetilde{a}_{ks}(y)\zeta ^{k}\zeta ^{s}\mbox{,}  \tag{5.2}
\end{equation}%
at $y=0$ for $\zeta \in 
%TCIMACRO{\U{211d} }%
%BeginExpansion
\mathbb{R}
%EndExpansion
^{n}$, where $\widetilde{\Gamma }_{jk}^{^{i}}$ are the Christoffel symbols
of the metric $g$ in the new coordinates $y$. We note that we will use the
symbols $\left( x,\xi \right) ,$ $u(x,\xi ),$ $\gamma (x,\xi ),$ $a_{ks}(x)$
instead of $\left( y,\zeta \right) ,$ $\widetilde{u}\left( y,\zeta \right) $%
, $\widetilde{\gamma }(y,\zeta )$, $\widetilde{a}_{ks}(y)$ in the
corresponding places of this paper. Since Lemmata 1-3 for the function $%
u(x,\xi )$ are also valid for the function $\widetilde{u}(y,\zeta )$,
without fear of confusion, we can change the symbols indicated above.

On the other hand, it is known that the equation of the geodesic in the
Euclidean metric has the form: $z\left( x_{0},\xi ,t\right) =x_{0}+\xi t,$ $%
(\xi \neq 0)$ and we readily see that the second derivatives of the
function\ $\overset{.}{z}\left( x_{0},\xi ,t\right) =\frac{d}{dt}z\left(
x_{0},\xi ,t\right) $ with respect to the parameter $\xi ^{s},$ $(1\leq
s\leq n)$ are bounded for $\left\vert \xi \right\vert \leq 1/2$ and $t\in
(0,+\infty ).$ At this point, it is worth to note that if the equation of
the geodesic is rewritten in the form $z\left( x_{0},\nu ,t\right)
=x_{0}+\nu t,$ where $\nu =\xi /\left\vert \xi \right\vert ,$ then the
second order derivatives of the function\ $\overset{.}{z}\left( x_{0},\nu
,t\right) =\xi /\left\vert \xi \right\vert $ with respect to the parameter $%
\xi ^{s}$ are unbounded for $\left\vert \xi \right\vert \leq 1/2$ and $t\in
(0,+\infty ).$ This unboundedness is connected with the introduction of the
new parameter $\nu =\xi /\left\vert \xi \right\vert $ and the unboundedness
of the second derivatives of the function $\left\vert \xi \right\vert $ for $%
\left\vert \xi \right\vert \leq 1/2.$ However in below, when we investigate
the boundedness of the derivatives with respect to the parameter $\xi ^{s}$
of the funcions $\overset{.}{z}\left( x_{0},\xi ,t\right) ,$ $\partial _{\xi
^{i}}\partial _{\xi ^{1}}u,$ $\partial _{x^{j}}\partial _{\xi ^{i}}\partial
_{\xi ^{1}}u$ in a neighbourhood of $\xi =0,1\leq i\leq n$, we need not pass
to the parameter $\nu =\xi /\left\vert \xi \right\vert $ and work only with
the parameter $\xi .$

\textbf{Lemma 4}. \textit{In the new coordinate system, the functions}%
\begin{equation*}
\partial _{\xi ^{i}}\partial _{\xi ^{1}}u(0,\xi )\text{\textit{,}}\mathit{\ }%
\partial _{\xi ^{i}}\partial _{\xi ^{1}}\partial _{x^{j}}u(0,\xi )\text{, }%
1\leq i,\text{ }j\leq n\text{,}
\end{equation*}%
\textit{\ are bounded on the set }$G$\textit{\ for }$\left\vert \xi
\right\vert \leq 1/2$. \textit{Moreover,} \textit{the functions }%
\begin{gather*}
u(0,\xi ^{1},\xi ^{\prime }),\text{ }\partial _{\xi ^{i}}u(0,\xi ^{1},\xi
^{\prime }),\text{ }\partial _{\xi ^{i}}\partial _{x^{j}}u(0,\xi ^{1},\xi
^{\prime }),\text{ }\partial _{\xi ^{i}}\partial _{\xi ^{1}}u(0,\xi ^{1},\xi
^{\prime }), \\
\partial _{x^{j}}\partial _{\xi ^{i}}\partial _{\xi ^{1}}u(0,\xi ^{1},\xi
^{\prime }),1\leq i,j\leq n,
\end{gather*}%
tend to zero in $L_{2}(\mathbb{R}_{\xi ^{1}}^{1})\ $as $\xi ^{\prime
}\rightarrow 0$.

\textbf{Proof.}

First, let us show that if up to the second order derivatives of the
functions $\overset{.}{z}^{i}\left( x,\xi ,t\right) $\ with respect to $\xi $
are bounded, then\ it is possible to prove that the functions $\partial
_{\xi ^{i}}\partial _{\xi ^{1}}u$ are bounded for $0<|\xi |\leq 1/2$ and $%
1\leq i\leq n$. \newline
Since supp $b\subset D$ and the solution of Problem (3.2)-(3.3) has property
(3.4), the integral in (4.2) can be considered on the finite interval $\left[
0,\frac{d_{0}}{\left\vert \xi \right\vert }\right] ,$ where $d_{0}$ is the
diameter of $D$ in the metric $g=\left( g_{ij}(x)\right) $. Hence we can
write (4.2) as 
\begin{equation}
I_{ij}(x,\xi )=\int_{0}^{\frac{d_{0}}{\left\vert \xi \right\vert }}b\left(
z\left( x,\xi ,t\right) \right) \overset{.}{z}^{i}\left( x,\xi ,t\right) 
\overset{.}{z}^{j}\left( x,\xi ,t\right) dt\text{.}  \tag{5.3}
\end{equation}%
Moreover, since the metric $g=\left( g_{ij}\right) $ is written down in the
semi-geodesic coordinates, the equality $\overset{.}{z}^{k}\left( x,\xi
^{1},0,t\right) =0$ ($\xi ^{1}\neq 0$) is satisfied for each $k$ $(2\leq
k\leq n)$ (see the proof of Lemma 1 in Appendix) and thus $\partial _{\xi
^{1}}z^{k}\left( x,\xi ^{1},0,t\right) =0.$ Then, for each fixed $x\in D$
and $0<|\xi ^{1}|\leq 1$, by Taylor's formula with respect to the variable $%
\xi ^{\prime }=(\xi ^{2},...,\xi ^{n})$ (e.g., \cite{Fikh}, p. 285), we
obtain the equalities 
\begin{equation}
\overset{.}{z}^{k}\left( x,\xi ,t\right) =\underset{i=2}{\overset{n}{\sum }}%
\xi ^{i}\partial _{\xi ^{i}}\overset{.}{z}^{k}\left( x,\xi ^{1},\xi ^{\prime
}\theta _{1},t\right) ,  \tag{5.4}
\end{equation}%
\begin{equation}
\partial _{\xi ^{1}}\overset{.}{z}^{k}\left( x,\xi ,t\right) =\underset{i=2}{%
\overset{n}{\sum }}\xi ^{i}\partial _{\xi ^{i}}\partial _{\xi ^{1}}\overset{.%
}{z}^{k}\left( x,\xi ^{1},\xi ^{\prime }\theta _{2},t\right) ,  \tag{5.5}
\end{equation}%
where $0<\theta _{m}\left( x,\xi ,t\right) <1,$ $m=1,2.$ If we assume the
boundedness of the functions $\partial _{\xi ^{i}}\overset{.}{z}^{k}$ and $%
\partial _{\xi ^{i}}\partial _{\xi ^{j}}\overset{.}{z}^{k}$, then it follows
from (5.4) and (5.5) that%
\begin{equation}
\left\vert \overset{.}{z}^{k}\left( x,\xi ,t\right) \right\vert \leq
M_{1}\left\vert \xi ^{\prime }\right\vert ,\mbox{ }\left\vert \partial _{\xi
^{1}}\overset{.}{z}^{k}\left( x,\xi ,t\right) \right\vert \leq
M_{2}\left\vert \xi ^{\prime }\right\vert  \tag{5.6}
\end{equation}%
for each $k\in \{2,...,n\}$, $x\in D,$ $t\in (0,+\infty )$ and for $%
\left\vert \xi \right\vert \leq 1/2$ ($\xi ^{1}\neq 0$), where $M_{1}$ and $%
M_{2}\geq 0$ do not depend on $\xi $ and $t$. Furthermore we assume the
boundedness of the functions $\partial _{\xi ^{i}}z^{s},$ $\partial _{\xi
^{j}}\partial _{\xi ^{i}}z^{s}$ $(1\leq s\leq n)$ and note that

(i) $b\in C^{5}(D),$ supp $b\subset D$,

(ii) Inequalities (5.6) hold,

(iii) $z\left( x,\xi ,t\right) $ is of $C^{5}$-class with respect to $t$ and 
$\xi \neq 0$,

(iv) Boundedness\ of the derivatives up to order 2 of the functions $\overset%
{.}{z}^{k}\left( x,\xi ,t\right) ,$\ $z^{s}\left( x,\xi ,t\right) $ $(2\leq
k\leq n,$ $1\leq s\leq n)$ with respect to $\xi ^{i}$ $(1\leq i\leq n)$\ for 
$|\xi |\leq 1/2$ with $\xi ^{1}\neq 0$,

(v) $\partial _{\xi ^{1}}z^{k}\left( x,\xi ^{1},0,t\right) =0$, $2\leq k\leq
n$. \newline
Then by (5.3) we see that the functions $\partial _{\xi ^{s}}\partial _{\xi
^{1}}I_{ij}$ are bounded for $\left\vert \xi \right\vert \leq 1/2,$ $1\leq
s\leq n$, $2\leq i,j\leq n$. On the other hand, the boundedness of the
functions for $\left\vert \xi \right\vert \geq 1/2$ ($\xi ^{\prime }\in
G^{\prime })$ follows from Lemma 1.

If we assume the boundedness of the functions $\overset{.}{\partial _{\xi
^{i}}\partial _{x^{j}}\overset{.}{z}^{k}},$ $\partial _{\xi ^{i}}\partial
_{\xi ^{s}}\partial _{x^{j}}\overset{.}{z}^{k},$ $\partial _{\xi
^{i}}\partial _{x^{j}}z^{s},$ $\partial _{\xi ^{l}}\partial _{\xi
^{i}}\partial _{x^{j}}z^{s}$ $(2\leq k\leq n,$ $1\leq i,$ $j,$ $s,$ $l\leq
n) $, then we can similarly prove the boundedness of the functions $\partial
_{\xi ^{i}}\partial _{\xi ^{1}}\partial _{x^{j}}I_{ij}$ for $\xi \in G$ $%
(1\leq i,j\leq n)$. As a result, since $u(x,\xi )$ is a finite sum of the
functions of the form $I_{ij}(x,\xi )$, the functions $\partial _{\xi
^{i}}\partial _{\xi ^{1}}u,$ $\partial _{\xi ^{i}}\partial _{\xi
^{1}}\partial _{x^{j}}u$ $(1\leq i,j\leq n)$ are bounded on the set $G$ for $%
x\in D$ and $|\xi |\leq 1/2$ with $\xi ^{1}\neq 0$, when the conditions for
the functions $z^{s},$ $\overset{.}{z}^{k}$ indicated above are satisfied.

In particular, since the derivatives of the functions $\overset{.}{z}%
^{k}\left( x,\xi ,t\right) $ with respect to $\xi $ up to order 2 are
bounded for the metric $g_{ij}=\delta _{ij}$, it is clear that the same
boundedness properties are valid for $\overset{.}{z}^{k}\left( 0,\xi
,t\right) $ in the Riemannian coordinates with center $\widetilde{x}_{0}\in
D $. Moreover, by the formula $y^{i}=\xi ^{i}t$ $(1\leq i\leq n)$ which is
satisfied in the new coordinate system, we have $z^{s}(0,\xi ,t)=z^{s}(0,\xi
,0)$ $\left( 1\leq s\leq n\right) $. Therefore, it is easy to see that the
functions $\partial _{\xi ^{i}}\partial _{x^{j}}z^{s}(0,\xi ,t),$ $\partial
_{\xi ^{l}}\partial _{\xi ^{i}}\partial _{x^{j}}z^{s}(0,\xi ,t),$ $\partial
_{\xi ^{i}}\partial _{x^{j}}\overset{.}{z}^{k}(0,\xi ,t)$, $\partial _{\xi
^{i}}\partial _{\xi ^{s}}\partial _{x^{j}}\overset{.}{z}^{k}(0,\xi ,t)$ $%
(2\leq k\leq n,$ $1\leq i,j,s,l\leq n)$ are bounded on the set $G$ for $%
\left\vert \xi \right\vert \leq 1/2.$ This implies that the functions $%
\partial _{\xi ^{i}}\partial _{\xi ^{1}}u(0,\xi )$, $\partial _{\xi
^{i}}\partial _{\xi ^{1}}\partial _{x^{j}}u(0,\xi )$ $(1\leq i,j\leq n) $
are bounded on the set $G$ for $\left\vert \xi \right\vert \leq 1/2.$

Now, let us prove the second assertion of the lemma. Since the unique
solution to Problem (3.2)-(3.3) for $\nu ^{0}=(1,0,...,0)\in 
%TCIMACRO{\U{211d} }%
%BeginExpansion
\mathbb{R}
%EndExpansion
^{n}$ is given by $z\left( x,\nu ^{0},t\right) =x+t\nu ^{0}$, we obtain that 
$\overset{.}{z}^{k}\left( x,\xi ^{1},0,t\right) =0$ for $\xi ^{\prime }=0$, $%
2\leq k\leq n$. Then by (3.1), we have $u(x,\xi ^{1},0)=\partial _{\xi
^{i}}u(x,\xi ^{1},0)=0$ $(1\leq i\leq n)$ and $\partial _{\xi ^{1}}\partial
_{\xi ^{i}}u(x,\xi ^{1},0)=\partial _{\xi ^{1}}\partial _{\xi ^{i}}\partial
_{x^{j}}u(x,\xi ^{1},0)=0$ $(1\leq i,j\leq n)$. Hence, taking into account
the last equalities, the first conclusion of the lemma which we proved above
and the analogous formula to (4.2), we see that the functions 
\begin{equation*}
u(0,\xi ^{1},\xi ^{\prime }),\,\partial _{\xi ^{i}}u(0,\xi ^{1},\xi ^{\prime
}),\,\partial _{\xi ^{i}}\partial _{x^{j}}u(0,\xi ^{1},\xi ^{\prime
}),\,\partial _{\xi ^{i}}\partial _{\xi ^{1}}u(0,\xi ^{1},\xi ^{\prime }),%
\text{ }\partial _{x^{j}}\partial _{\xi ^{i}}\partial _{\xi ^{1}}u(0,\xi
^{1},\xi ^{\prime })
\end{equation*}%
tend to zero in $L_{2}(\mathbb{R}_{\xi ^{1}}^{1})\ $as $\xi ^{\prime
}\rightarrow 0.$

\textbf{Lemma 5. }\textit{Let the conditions of Lemmata 1 and 3 be
satisfied. \textit{We set }}$\widehat{u}(0,\eta ,\xi ^{\prime })=p(0,\eta
,\xi ^{\prime })+\sqrt{-1}q(0,\eta ,\xi ^{\prime })$\textit{\textit{.} Then
the functions }$p(0,\eta ,\xi ^{\prime }),$ $\partial _{x^{j}}p(0,\eta ,\xi
^{\prime }),$\textit{\ }$\partial _{\xi ^{k}}p(0,\eta ,\xi ^{\prime }),$%
\textit{\ }$\partial _{x^{j}}\partial _{\xi ^{k}}p(0,\eta ,\xi ^{\prime }),$ 
$\eta q(0,\eta ,\xi ^{\prime })$\textit{\ tend to zero as }$\xi ^{\prime
}\rightarrow 0$ \mbox{in} $L_{1}(\mathbb{R}_{\eta }^{1})$ \textit{for }$%
2\leq k\leq n$\textit{\ and }$1\leq j\leq n$\textit{.}

\textbf{Proof. }Since the Fourier transform is continuous in $L_{2}$($%
%TCIMACRO{\U{211d} }%
%BeginExpansion
\mathbb{R}
%EndExpansion
_{\xi ^{1}}^{1}$), by Lemma 4 the functions $\widehat{u}(0,\eta ,\xi
^{\prime })$, $\partial _{\xi ^{k}}\widehat{u}(0,\eta ,\xi ^{\prime })$, $%
\eta \widehat{u}(0,\eta ,\xi ^{\prime })$, $\eta \partial _{\xi ^{k}}%
\widehat{u}(0,\eta ,\xi ^{\prime })$, $\eta \partial _{x^{j}}\widehat{u}%
(0,\eta ,\xi ^{\prime })$, $\partial _{\xi ^{k}}\partial _{x^{j}}\widehat{u}%
(0,\eta ,\xi ^{\prime })$, $\eta ^{2}\partial _{x^{j}}\widehat{u}(0,\eta
,\xi ^{\prime })$, $\eta \partial _{\xi ^{k}}\partial _{x^{j}}\widehat{u}%
(0,\eta ,\xi ^{\prime })$, $2\leq k\leq n$, $1\leq j\leq n$, tend to zero in 
$L_{2}$($%
%TCIMACRO{\U{211d} }%
%BeginExpansion
\mathbb{R}
%EndExpansion
_{\eta }^{1}$) as $\xi ^{\prime }\rightarrow 0$. Then, since 
\begin{equation*}
\int_{1}^{\infty }\left\vert \partial _{\xi ^{k}}\partial
_{x^{j}}p\right\vert d\eta \leq \left( \int_{1}^{\infty }\eta ^{2}\left(
\partial _{\xi ^{k}}\partial _{x^{j}}p\right) ^{2}d\eta \right) ^{\frac{1}{2}%
}\left( \int_{1}^{\infty }\frac{1}{\eta ^{2}}d\eta \right) ^{\frac{1}{2}},
\end{equation*}%
\begin{equation*}
\int_{0}^{1}\left\vert \partial _{\xi ^{k}}\partial _{x^{j}}p\right\vert
d\eta \leq \left( \int_{0}^{1}\left( \partial _{\xi ^{k}}\partial
_{x^{j}}p\right) ^{2}\mbox{ }d\eta \right) ^{\frac{1}{2}}
\end{equation*}%
and%
\begin{equation*}
\int_{1}^{\infty }\eta ^{2}\left( \partial _{\xi ^{k}}\partial
_{x^{j}}p\right) ^{2}d\eta \rightarrow 0,\text{ }\int_{0}^{1}\left( \partial
_{\xi ^{k}}\partial _{x^{j}}p\right) ^{2}d\eta \rightarrow 0,
\end{equation*}%
we see that the functions $\partial _{\xi ^{k}}\partial _{x^{j}}p(0,\eta
,\xi ^{\prime })$ tend to zero in $L_{1}$($0,\infty $) as $\xi ^{\prime
}\rightarrow 0$. By the same argument, it can be proved that the function $%
\partial _{\xi ^{k}}\partial _{x^{j}}p(0,\eta ,\xi ^{\prime })$ tends to
zero in $L_{1}(-\infty ,0)$ as $\xi ^{\prime }\rightarrow 0$.

Similarly one can prove that the functions $p(0,\eta,\xi^{\prime })$, $%
\partial_{x^{j}}p(0,\eta,\xi^{\prime })$, $\partial_{\xi
^{k}}p(0,\eta,\xi^{\prime })$, $\eta q(0,\eta,\xi^{\prime })$ tend to zero
in $L_{1}(%
%TCIMACRO{\U{211d} }%
%BeginExpansion
\mathbb{R}
%EndExpansion
_{\eta }^{1})$ as $\xi ^{\prime }\rightarrow 0.$

\section{\textbf{The proof of\ the main result}}

In this section, we prove Theorem 1. As it was noted in Section 3, the
uniqueness of the solution to Problem 1 follows from the uniqueness of the
solution to Problem 2 in the class $C_{0}^{5}(D)$. Therefore, we consider
Problem 2 with homogeneous boundary data below.

\textbf{Proof of Theorem 1.}

First we recall that$\ \widehat{u}=p + \sqrt{-1}q$ and $\Gamma
_{jk}^{i}(0)=0 $.\ Then by (4.3) and (4.6), in the new coordinates, we have%
\begin{equation}
\partial _{\eta }\partial _{x^{1}}q=-2\pi \delta (\eta )\underset{k,j=2}{%
\overset{n}{\sum }}a_{jk}\left( 0\right) \xi ^{^{k}}\xi ^{j}+\digamma _{2}, 
\tag{6.1}
\end{equation}%
for $x=0$, where $\digamma _{2}=\underset{j=2}{\overset{n}{%
\mathop{\textstyle \sum }}}\xi ^{j}\partial _{x^{j}}p$. By (6.1), we obtain%
\begin{equation}
\partial _{\eta }\left( \partial _{x^{1}}q(0,\eta ,\xi _{i}^{\prime
})\right) =-2\pi \delta \left( \eta \right) a_{ii}(0)\varepsilon
^{2}+F_{2i}(0,\eta ,\xi _{i}^{\prime })\text{,}  \tag{6.2}
\end{equation}%
for $\xi ^{\prime }=\xi _{i}^{\prime }\in G^{\prime }$, $2\leq i\leq n$,
where $\xi _{i}^{\prime }=\varepsilon \xi _{i}^{\prime }(1),$ $\xi
_{i}^{\prime }(1)=(\underset{i-2}{\underbrace{0,\cdots ,0,}}1,0,\cdots
,0)\in 
%TCIMACRO{\U{211d} }%
%BeginExpansion
\mathbb{R}
%EndExpansion
^{n-1},$ $\varepsilon >0,$ $\xi ^{\prime }=\left( \xi ^{2},\xi ^{3},\cdots
,\xi ^{n}\right) $ and\ $F_{2i}(0,\eta ,\xi _{i}^{\prime })=\varepsilon
\partial _{x^{i}}p(0,\eta ,\xi _{i}^{\prime })$.

On the other hand, we know (e.g., Theorem 3.1.3 on p. 56 in H\"{o}rmander 
\cite{Horm})

\textbf{Lemma 6. }\textit{Let a function }$U(y)$ \textit{defined in an open
set }$Y\subset \mathbb{R}$, \textit{\ belong to space }$C^{1}(Y/\left\{
y_{0}\right\} )$\textit{\ for some }$y_{0}\in Y$\textit{\ and let a function 
}$V(y)$ \textit{coincide with }$\dfrac{dU\left( y\right) }{dy}$\textit{\ for}
$y\neq y_{0}$ \textit{and be integrable on some neighbourhood of }$y_{0}$. 
\textit{Then the limits}%
\begin{equation*}
U(y_{0}\pm 0):=\underset{y\rightarrow y_{0}\pm }{\lim }U(y)
\end{equation*}%
\textit{exist and}%
\begin{equation*}
\dfrac{dU\left( y\right) }{dy}=V(y)+(U(y_{0}+0)-U(y_{0}-0))\delta \left(
y_{0}\right) \text{.}
\end{equation*}%
Lemma 3\ shows that (taking into account Remark 2) for fixed $(0,\xi
^{\prime })\in D\times G^{\prime }$, the functions $U=\partial _{x^{1}}q$ and%
$\ V=\digamma _{2}$ satisfy the conditions of Lemma 6 for $y_{0}=0$, in
where the variable $y$ is replaced by $\eta $. Therefore, by (6.2) and Lemma
6, we conclude that%
\begin{equation}
U_{+q}(0,\xi _{i}^{\prime })-U_{-q}(0,\xi _{i}^{\prime })=-2\pi
a_{ii}(0)\varepsilon ^{2},  \tag{6.3}
\end{equation}%
where $U_{\pm q}(0,\xi _{i}^{\prime })=\partial _{x^{1}}q\left( 0,\pm 0,\xi
_{i}^{\prime }\right) $.

On the other hand, from (6.2) by Lemmata 2 and 3, it is not difficult to
obtain 
\begin{equation*}
U_{+q}(0,\xi _{i}^{\prime })=-\underset{0}{\overset{\infty }{\int }}%
F_{2i}(0,\eta ,\xi _{i}^{\prime })d\eta ,\ U_{-q}(0,\xi _{i}^{\prime })=%
\underset{-\infty }{\overset{0}{\int }}F_{2i}(0,\eta ,\xi _{i}^{\prime
})d\eta \text{.}
\end{equation*}%
Consequently, we have%
\begin{equation}
U_{+q}(0,\xi _{i}^{\prime })-U_{-q}(0,\xi _{i}^{\prime })=-\varepsilon 
\underset{-\infty }{\overset{\infty }{\int }}\partial _{x^{_{i}}}p(0,\eta
,\xi _{i}^{\prime })d\eta \text{.}  \tag{6.4}
\end{equation}%
\ By Lemma 5,%
\begin{equation*}
\underset{-\infty }{\overset{\infty }{\int }}\partial _{x^{s}}p(0,\eta
,0)d\eta =0\text{, }1\leq s\leq n\text{,}
\end{equation*}%
and Lemmata 2, 3 and 5 yield 
\begin{equation*}
\partial _{\xi ^{i}}\partial _{x^{_{i}}}p\in L_{1}\left( 
%TCIMACRO{\U{211d} }%
%BeginExpansion
\mathbb{R}
%EndExpansion
_{\eta }^{1}\right) \cap L_{2}\left( 
%TCIMACRO{\U{211d} }%
%BeginExpansion
\mathbb{R}
%EndExpansion
_{\eta }^{1}\right) \cap C\left( \Delta _{\eta }^{\rho }\times G^{\prime
}\right) ,\text{ }\rho =-1,1\text{.}
\end{equation*}%
Hence, for $x=0,$ the mean value theorem on interval $\left[ 0,\varepsilon %
\right] $ (e.g., \cite{Fikh}, p. 186) implies that%
\begin{equation}
\underset{-\infty }{\overset{\infty }{\int }}\partial _{x^{i}}p(0,\eta ,\xi
_{i}^{\prime })d\eta =\varepsilon \underset{-\infty }{\overset{\infty }{\int 
}}\partial _{\xi ^{i}}\partial _{x^{_{i}}}p(0,\eta ,\xi _{i}^{\prime }\theta
_{1})d\eta ,  \tag{6.5}
\end{equation}%
where $\theta _{1}\in (0,1)$ is a constant depending on $\xi _{i}^{\prime }$%
. By using (6.5) in (6.4), we obtain 
\begin{equation}
U_{+q}(0,\xi _{i}^{\prime })-U_{-q}(0,\xi _{i}^{\prime })=q_{i}(0,\xi
_{i}^{\prime })\varepsilon ^{2},  \tag{6.6}
\end{equation}%
where%
\begin{equation*}
q_{i}(0,\xi _{i}^{\prime })=\underset{-\infty }{-\overset{\infty }{\int }}%
\partial _{\xi ^{i}}\partial _{x^{_{i}}}p(0,\eta ,\xi _{i}^{\prime }\theta
_{1})d\eta .
\end{equation*}%
Equalities (6.3) and (6.6) show that $-2\pi a_{ii}(0)=q_{i}(0,\xi
_{i}^{\prime }),$ $(2\leq i\leq n).$ Then by Lemma 5, we have $q_{i}(0,\xi
_{i}^{\prime })\rightarrow 0$ as $\xi _{i}^{\prime }\rightarrow 0$ and thus $%
a_{ii}(0)=0,$ $(2\leq i\leq n).$

In order to complete the proof, let $\xi _{ij}^{\prime }(1)=(\underset{i-2}{%
\underbrace{0,\cdots ,0}},1,\underset{j-i-1}{\underbrace{0,\cdots ,0},}1,%
\underset{n-j}{\underbrace{0,\cdots ,0}})\in 
%TCIMACRO{\U{211d} }%
%BeginExpansion
\mathbb{R}
%EndExpansion
^{n-1},$ $i\neq j,$ and $\xi _{ij}^{\prime }=\varepsilon \xi _{ij}^{\prime
}(1)\in G^{\prime }$.\ In this case, by $a_{ii}\left( 0\right) =0$, $%
a_{ij}\left( 0\right) =a_{ji}\left( 0\right) $ and (6.1), we have 
\begin{equation}
\partial _{\eta }\left( \partial _{x^{1}}q\left( 0,\eta ,\xi _{ij}^{\prime
}\right) \right) =-4\pi \delta \left( \eta \right) a_{ij}(0)\varepsilon
^{2}+F_{2ij}\left( 0,\eta ,\xi _{ij}^{\prime }\right)  \tag{6.7}
\end{equation}%
for $\xi ^{\prime }=\xi _{ij}^{\prime }$, where 
\begin{equation*}
F_{2ij}(0,\eta ,\xi _{ij}^{\prime })=\varepsilon (\partial _{x^{i}}p(0,\eta
,\xi _{ij}^{\prime })+\partial _{x^{j}}p(0,\eta ,\xi _{ij}^{\prime }))
\end{equation*}%
and the indices $i,j\in \{2,...,n\}$ are fixed. Then, in the same way as
above, equation (6.7) and Lemma 6 yield%
\begin{equation}
U_{+q}(0,\xi _{ij}^{\prime })-U_{-q}(0,\xi _{ij}^{\prime })=-4\pi
a_{ij}(0)\varepsilon ^{2},  \tag{6.8}
\end{equation}%
where $U_{\pm q}(0,\xi _{ij}^{\prime })=\partial _{x^{1}}q\left( 0,\pm 0,\xi
_{ij}^{\prime }\right) $. From (6.7),\ by virtue of Lemmata 2, 3 we have%
\begin{equation*}
U_{+q}(0,\xi _{ij}^{\prime })=-\underset{0}{\overset{\infty }{\int }}%
F_{2ij}\left( 0,\eta ,\xi _{ij}^{\prime }\right) d\eta ,\mbox{ \ }%
U_{-q}(0,\xi _{ij}^{\prime })=\underset{-\infty }{\overset{0}{\int }}%
F_{2ij}\left( 0,\eta ,\xi _{ij}^{\prime }\right) d\eta \text{.}
\end{equation*}%
Hence,%
\begin{equation}
U_{+q}(0,\xi _{ij}^{\prime })-U_{-q}(0,\xi _{ij}^{\prime })=-\varepsilon 
\underset{-\infty }{\overset{\infty }{\int }}(\partial _{x^{i}}p(0,\eta ,\xi
_{ij}^{\prime })+\partial _{x^{j}}p(0,\eta ,\xi _{ij}^{\prime }))d\eta . 
\tag{6.9}
\end{equation}%
Recalling that (see Lemma 5)%
\begin{equation*}
\underset{-\infty }{\overset{\infty }{\int }}\partial _{x^{m}}p(0,\eta
,0)d\eta =0
\end{equation*}%
$(1\leq m\leq n),$\ by the mean value theorem on interval $\left[
0,\varepsilon \right] $ (e.g., \cite{Fikh}, p. 186), we obtain 
\begin{equation*}
\underset{-\infty }{\overset{\infty }{\int }}\partial _{x^{s}}p(0,\eta ,\xi
_{ij}^{\prime })d\eta =\varepsilon \underset{-\infty }{\overset{\infty }{%
\int }}(\partial _{\xi ^{i}}\partial _{x^{s}}p(0,\eta ,\xi _{ij}^{\prime
}\theta _{1}^{s})+\partial _{\xi ^{j}}\partial _{x^{s}}p(0,\eta ,\xi
_{ij}^{\prime }\theta _{1}^{s}))d\eta ,
\end{equation*}%
where\ $s=i,j;$ and $0<\theta _{1}^{s}\left( \xi _{ij}^{\prime }\right) <1.$
Then, by the last equality and (6.9), we have%
\begin{equation}
U_{+q}(0,\xi _{ij}^{\prime })-U_{-q}(0,\xi _{ij}^{\prime })=q_{ij}(0,\xi
_{ij}^{\prime })\varepsilon ^{2},  \tag{6.10}
\end{equation}%
where 
\begin{equation*}
q_{ij}(0,\xi _{ij}^{\prime })=\underset{-\infty }{-\overset{\infty }{\int }}%
\left( \partial _{\xi ^{i}}\partial _{x^{_{i}}}p(0,\eta ,\xi _{ij}^{\prime
}\theta _{1}^{s})+2\partial _{\xi ^{j}}\partial _{x^{_{i}}}p(0,\eta ,\xi
_{ij}^{\prime }\theta _{1}^{s})+\partial _{\xi ^{j}}\partial
_{x^{_{j}}}p(0,\eta ,\xi _{ij}^{\prime }\theta _{1}^{s})\right) d\eta .
\end{equation*}%
Equalities (6.8) and (6.10) imply $-4\pi a_{ij}\left( 0\right) =q_{ij}\left(
0,\xi _{ij}^{\prime }\right) $. Then, by an argument similar to the proof of
the relation $q_{i}(0,\xi _{i}^{\prime })\rightarrow 0$ as $\xi _{i}^{\prime
}\rightarrow 0$, we can prove that $q_{ij}(0,\xi _{ij}^{\prime })\rightarrow
0$ as $\xi _{i}^{\prime }\rightarrow 0,$ therefore, $a_{ij}\left( 0\right)
=0,$ $(2\leq i,$ $j\leq n)$.

In fact, keeping in mind the change of variables $x\mapsto y$ (see the
beginning of Section 5), instead of $a_{ks}(0)$, it suffices to consider $%
\widetilde{a}_{ks}(0)$, and we have $a_{ks}(0)=\widetilde{a}_{ks}(0)=a_{ij}(%
\widetilde{x}_{0})\left. \partial _{y_{k}}x^{i}\partial
_{y_{s}}x^{j}\right\vert _{y=0},$ $2\leq k,$ $s\leq n.$

Then from $\widetilde{a}_{ks}(0)=0$ $(2\leq k,$ $s\leq n)$, it follows that $%
\widetilde{a}_{ks}(0)=A_{is}\left. \partial _{y^{k}}x^{i}\right\vert
_{y=0}=0 $ for each $2\leq s\leq n,$ where $A_{is}=a_{ij}(\widetilde{x}%
_{0})\left. \partial _{y_{s}}x^{j}\right\vert _{y=0}.$ Since the Jacobian
det $(\partial _{y_{s}}x^{j})\neq 0$, we have $A_{is}=0$ for $2\leq i,$ $%
s\leq n$ and so $a_{ij}(\widetilde{x}_{0})=0$ for $2\leq i,$ $j\leq n.$
Since $\widetilde{x}_{0}\in D$ is arbitrary, we obtain $a_{ij}=0$ for $2\leq
i,j\leq n$ in the domain $D.$ Hence the proof of theorem 1 is completed.

\section{Appendix}

In this section, we prove Lemmata 1-3. For the sake of simplicity, we set 
\begin{equation*}
I(x,\xi ):=I_{ij}(x,\xi )
\end{equation*}%
in (4.1) after fixing the indices $i,j\in \left\{ 2,...,n\right\} .$

\textbf{Proof of Lemma 1. }Similarly to (5.3), we recall that the last
integral in (4.2) is considered on the finite interval $\left[ 0,d_{0}\right]
$. Due to the condition\ $g_{ij}\in C^{6}(D)$, it follows from the theory of
ordinary differential equations that the solution $z(x,\nu ,t)$ of Problem
(3.2)-(3.3) belongs to the space $C^{5}(\Omega (d_{0}))$, where we recall
that $\Omega \left( d_{0}\right) =\left\{ \left( x,\nu ,t\right) |x\in D,\nu
\in S^{n}(x),t\in \left[ 0,d_{0}\right] \right\} $ and $S^{n}(x)$ is the
unit sphere with the metric $g$ around $x\in D.$ Hence, taking into account
the conditions $b\in C^{5}\left( 
%TCIMACRO{\U{211d} }%
%BeginExpansion
\mathbb{R}
%EndExpansion
^{n}\right) ,$ $\xi \in G,$ $\xi ^{\prime }\neq 0$ and equality (4.2), it
follows that $\partial _{\xi }^{\beta }I,$ $\partial _{\xi }^{\beta
}\partial _{x^{s}}I\in C(\Omega )$ for $0\leq \left\vert \beta \right\vert
\leq 4,$ $1\leq s\leq n.$

In order to prove the second assertion of Lemma 1, we investigate the
behavior of $\left\vert \xi \right\vert \overset{.}{z}^{k}\left( x,\nu ,\tau
\right) $ and its derivatives with respect to $\xi ^{j}$ when\ $\xi
^{1}\longrightarrow \infty ,$ where $1\leq j\leq n$ and $2\leq k\leq n$.

Let $\xi ^{1}=1/\mu $. Then the vector$\ \nu =\xi /\left\vert \xi
\right\vert \in S^{n}(x)$ tends to $\nu ^{0}=\left( 1,0,...,0\right) \in 
%TCIMACRO{\U{211d} }%
%BeginExpansion
\mathbb{R}
%EndExpansion
^{n}$ as$\ \xi ^{1}\longrightarrow +\infty ,$ (i.e. as$\ \mu \longrightarrow
+0$). Therefore, as it is known from the theory of ordinary differential
equations (\cite{Pet}), the unique solution of Problem (3.2)-(3.3) tends to
the solution $z\left( x,\nu ^{0},t\right) $ in $C^{2}[0,d_{0}]$ as\ $\mu
\longrightarrow +0.$ Since the metric $g=\left( g_{ij}\right) $ is given in
the semi-geodesic coordinates (then $\Gamma _{1s}^{1}=\Gamma _{11}^{s}=0,$ $%
1\leq s\leq n$) and the solution of Problem (3.2)-(3.3) is unique, we have $%
z\left( x,\nu ^{0},t\right) =\left( z^{1}\left( x,\nu ^{0},t\right)
,...,z^{n}\left( x,\nu ^{0},t\right) \right) ,$ where $z^{1}\left( x,\nu
^{0},t\right) =x^{1}+t,~z^{k}\left( x,\nu ^{0},t\right) =x^{k},$ $\overset{.}%
{z}\left( x,\nu ^{0},t\right) =(1,0,\cdots ,0),$ $\overset{.}{z}^{1}\left(
x,\nu ^{0},t\right) =1,$ $\overset{.}{z}^{k}\left( x,\nu ^{0},t\right) =0,$ $%
2\leq k\leq n$.

For $\xi ^{1}=1/\mu>0 $, we have%
\begin{equation}
\left\vert \xi \right\vert \overset{.}{z}^{k}\left( x,\frac{\xi ^{1}}{%
\left\vert \xi \right\vert },\cdots ,\frac{\xi ^{n}}{\left\vert \xi
\right\vert },t\right) =\frac{1}{\mu }\left\vert \xi ^{\prime }\right\vert
_{\mu }\overset{.}{z}^{k}\left( x,\frac{1}{\left\vert \xi ^{\prime
}\right\vert _{\mu }},\frac{\mu \xi ^{2}}{\left\vert \xi ^{\prime
}\right\vert _{\mu }},\cdots ,\frac{\mu \xi ^{n}}{\left\vert \xi ^{\prime
}\right\vert _{\mu }},t\right),  \tag{A.1}
\end{equation}%
where%
\begin{equation*}
\left\vert \xi ^{\prime }\right\vert _{\mu }=\left( 1+\mu ^{2}\underset{2}{%
\overset{n}{\sum }}g_{ij}\xi ^{i}\xi ^{j}\right) ^{1/2}\text{.}
\end{equation*}%
Hence, if $x\in D,$ $(\xi ^{2},...,\xi ^{n})\in G^{\prime },$ $t\in \left[
0,d_{0}\right] $ are fixed, then applying the mean value theorem (e.g., \cite%
{Fikh}, p.186) to the function $\overset{.}{z}^{k}\left( x,\frac{1}{%
\left\vert \xi ^{\prime }\right\vert _{\theta }},\frac{\theta \xi ^{2}}{%
\left\vert \xi ^{\prime }\right\vert _{\theta }},\cdots ,\frac{\theta \xi
^{n}}{\left\vert \xi ^{\prime }\right\vert _{\theta }},t\right) $ with
respect to $\theta $\ on the interval $\left[ 0,\mu \right] ,$ from the
equality $\overset{.}{z}^{k}\left( x,\nu ^{0},t\right) =0$ for $2\le k\le n$%
, we have 
\begin{equation}
\overset{.}{z}^{k}\left( x,\frac{1}{\left\vert \xi ^{\prime }\right\vert
_{\mu }},\frac{\mu \xi ^{2}}{\left\vert \xi ^{\prime }\right\vert _{\mu }}%
,\cdots ,\frac{\mu \xi ^{n}}{\left\vert \xi ^{\prime }\right\vert _{\mu }}%
,t\right) =\mu \partial _{\mu _{0}}\overset{.}{z}^{k}, \thinspace \thinspace
0<\mu_{0}<\mu \leq 1, \thinspace 2 \le k \le n,  \tag{A.2}
\end{equation}%
where $\partial _{\mu _{0}}\overset{.}{z}^{k}$ is the derivative of the
function $\overset{.}{z}^{k}\left( x,\frac{1}{\left\vert \xi ^{\prime
}\right\vert _{\theta }},\frac{\theta \xi ^{2}}{\left\vert \xi ^{\prime
}\right\vert _{\theta }},\cdots ,\frac{\theta \xi ^{n}}{\left\vert \xi
^{\prime }\right\vert _{\theta }},t\right) $ with respect to $\theta $\ at a
point $\theta =\mu _{0}.$ It is worth to note here that $\overset{.}{z}%
^{1}\left( x,\nu ^{0},t\right) =1$ and equality (A.2) is not valid for $k=1.$

Since $\overset{.}{z}^{k}\left( x,\nu ,t\right) \in C^{5}(\Omega (d_{0}))$
and $\left( \frac{1}{\left\vert \xi ^{\prime }\right\vert _{\mu }},\frac{\mu
\xi ^{2}}{\left\vert \xi ^{\prime }\right\vert _{\mu }},\cdots ,\frac{\mu
\xi ^{n}}{\left\vert \xi ^{\prime }\right\vert _{\mu }}\right) \in S^{n}(x),$
the function $\partial _{\mu }\overset{.}{z}^{k}\left( x,\frac{1}{\left\vert
\xi ^{\prime }\right\vert _{\mu }},\frac{\mu \xi ^{2}}{\left\vert \xi
^{\prime }\right\vert _{\mu }},\cdots ,\frac{\mu \xi ^{n}}{\left\vert \xi
^{\prime }\right\vert _{\mu }},t\right) $ is bounded on $\Omega \left(
d_{0}\right) $. Here we note that $\Omega \left( d_{0}\right) $ is closed
and bounded. Therefore,\ by (A.1) and (A.2), since the vector $\nu =\xi
/\left\vert \xi \right\vert \in S^{n}(x)$ tends to $\nu ^{0}=\left(
1,0,...,0\right) \in 
%TCIMACRO{\U{211d} }%
%BeginExpansion
\mathbb{R}
%EndExpansion
^{n}$ as\ $\xi ^{1}\longrightarrow +\infty $,$\ $we have%
\begin{equation}
\left\vert \mbox{ }\left\vert \xi \right\vert \mbox{ }\overset{.}{z}%
^{k}\left( x,\nu ,t\right) \right\vert \leq K_{1},  \tag{A.3}
\end{equation}%
for $2\leq k\leq n$ in the set $\Omega ,$ where $K_{1}>0$ is independent of $%
(x,\xi )\in D\times G$, but depends on the norm of the vector function $%
\overset{.}{z}\left( x,\nu ,t\right) $ in $C^{1}\left( \Omega \left(
d_{0}\right) \right) $ and the diameter of $G^{\prime }$. In the same way as
above, we can prove the last inequality for the case $\xi
^{1}\longrightarrow -\infty .$

It is not difficult to verify the following equalities%
\begin{eqnarray}
\partial _{\xi ^{1}}\left( \left\vert \xi \right\vert \overset{.}{z}%
^{k}\left( x,\frac{\xi }{\left\vert \xi \right\vert },t\right) \right) &=&%
\frac{\xi ^{1}}{\left\vert \xi \right\vert ^{2}}\left( \left\vert \xi
\right\vert \overset{.}{z}^{k}\left( x,\frac{\xi }{\left\vert \xi
\right\vert },t\right) \right)  \notag \\
&&+\left\vert \xi \right\vert \left( -\overset{n}{\underset{j=1}{\sum }}%
\partial _{\nu ^{j}}\overset{.}{z}^{k}\frac{\xi ^{j}\xi ^{1}}{\left\vert \xi
\right\vert ^{3}}+\frac{1}{\left\vert \xi \right\vert }\partial _{\nu ^{1}}%
\overset{.}{z}^{k}\right) ,  \TCItag{A.4}
\end{eqnarray}%
\begin{eqnarray}
&&\partial _{\xi ^{i}}\left( \left\vert \xi \right\vert \overset{.}{z}%
^{k}\left( x,\frac{\xi }{\left\vert \xi \right\vert },t\right) \right) =%
\frac{1}{\left\vert \xi \right\vert ^{2}}\underset{j=2}{\overset{n}{\sum }}%
g_{ij}\xi ^{j}\left( \left\vert \xi \right\vert \overset{.}{z}^{k}\left( x,%
\frac{\xi }{\left\vert \xi \right\vert },t\right) \right)  \notag \\
&&+\left\vert \xi \right\vert \left( -\overset{n}{\underset{s=1}{\sum }}%
\left( \partial _{\nu ^{s}}\overset{.}{z}^{k}\right) \frac{1}{\left\vert \xi
\right\vert ^{3}}\xi ^{s}\underset{j=2}{\overset{n}{\sum }}g_{ij}\xi ^{j}+%
\frac{1}{\left\vert \xi \right\vert }\partial _{\nu ^{i}}\overset{.}{z}%
^{k}\right) ,  \TCItag{A.5}
\end{eqnarray}%
\begin{eqnarray}
\partial _{\xi ^{1}}^{2}\left( \left\vert \xi \right\vert \overset{.}{z}%
^{k}\left( x,\frac{\xi }{\left\vert \xi \right\vert },t\right) \right) &=&%
\frac{1}{\left\vert \xi \right\vert }\overset{.}{z}^{k}-\frac{(\xi ^{1})^{2}%
}{\left\vert \xi \right\vert ^{3}}\overset{.}{z}^{k}-\overset{n}{\underset{%
j=2}{\sum }}\frac{\xi ^{j}}{\left\vert \xi \right\vert ^{2}}\partial _{\nu
^{j}}\overset{.}{z}^{k}-2\frac{\xi ^{1}}{\left\vert \xi \right\vert ^{2}}%
\partial _{\nu ^{1}}\overset{.}{z}^{k}  \notag \\
&&+2\xi ^{1}\overset{n}{\underset{j=1}{\sum }}\frac{\xi ^{j}\xi ^{1}}{%
\left\vert \xi \right\vert ^{4}}\partial _{\nu ^{j}}\overset{.}{z}^{k}+\frac{%
\xi ^{1}}{\left\vert \xi \right\vert ^{2}}\left( \overset{n}{\underset{j=1}{%
\sum }}\frac{-\xi ^{j}\xi ^{1}}{\left\vert \xi \right\vert ^{2}}\partial
_{\nu ^{j}}\overset{.}{z}^{k}+\partial _{\nu ^{1}}z^{k}\right)  \notag \\
&&\overset{n}{\underset{j=1}{+\sum }}\frac{\xi ^{j}\xi ^{1}}{\left\vert \xi
\right\vert ^{3}}\left( \sum\limits_{i=1}^{n}\frac{\xi ^{i}\xi ^{1}}{%
\left\vert \xi \right\vert ^{2}}\partial _{\nu ^{i}}\partial _{\nu ^{j}}%
\overset{.}{z}^{k}-\partial _{\nu ^{j}}\partial _{\nu ^{1}}\overset{.}{z}%
^{k}\right)  \notag \\
&&-\overset{n}{\underset{j=1}{\sum }}\frac{\xi ^{j}\xi ^{1}}{\left\vert \xi
\right\vert ^{3}}\partial _{\nu ^{1}}\partial _{\nu ^{j}}\overset{.}{z}^{k}+%
\frac{1}{\left\vert \xi \right\vert }\partial _{\nu ^{1}}^{2}\overset{.}{z}%
^{k},  \TCItag{A.6}
\end{eqnarray}%
\begin{eqnarray}
\partial _{\xi ^{i}}\partial _{\xi ^{1}}\left( \left\vert \xi \right\vert 
\overset{.}{z}^{k}\left( x,\frac{\xi }{\left\vert \xi \right\vert },t\right)
\right) &=&-\frac{\xi ^{1}}{\left\vert \xi \right\vert ^{3}}(\underset{j=2}{%
\overset{n}{\sum }}g_{ij}\xi ^{j})\overset{.}{z}^{k}-\left( \overset{n}{%
\underset{j=1}{\sum }}\overset{.}{z}_{\nu ^{j}}^{k}\frac{\xi ^{j}\xi ^{1}}{%
\left\vert \xi \right\vert ^{4}}\right) \left( \underset{j=2}{\overset{n}{%
\sum }}g_{ij}\xi ^{j}\right)  \notag \\
&&+\overset{.}{z}_{\nu ^{i}}^{k}\frac{\xi ^{1}}{\left\vert \xi \right\vert
^{2}}+2\left( \overset{n}{\underset{j=1}{\sum }}\overset{.}{z}_{\nu
^{j}}^{k}\xi ^{j}\xi ^{1}\right) \left( \frac{1}{\left\vert \xi \right\vert
^{4}}\underset{j=2}{\overset{n}{\sum }}g_{ij}\xi ^{j}\right)  \notag \\
&&-\overset{.}{z}_{\nu ^{i}}^{k}\frac{\xi ^{1}}{\left\vert \xi \right\vert
^{2}}+\overset{n}{\underset{j=1}{\sum }}\frac{\xi ^{j}\xi ^{1}}{\left\vert
\xi \right\vert ^{3}}\left( \overset{n}{\underset{m=1}{\sum }}\overset{.}{z}%
_{\nu ^{j}\nu ^{m}}^{k}\frac{1}{\left\vert \xi \right\vert ^{2}}\xi ^{m}%
\underset{s=2}{\overset{n}{\sum }}g_{is}\xi ^{s}-\overset{.}{z}_{\nu ^{i}\nu
^{j}}^{k}\right)  \notag \\
&&-\overset{n}{\underset{j=1}{\mathop{\textstyle \sum }}}\overset{.}{z}_{\nu
^{j}\nu ^{1}}^{k}\frac{1}{\left\vert \xi \right\vert ^{3}}\xi ^{j}\underset{%
j=2}{\overset{n}{\sum }}g_{ij}\xi ^{j}+\frac{1}{\left\vert \xi \right\vert }%
\overset{.}{z}_{\nu ^{i}\nu ^{1}}^{k},  \TCItag{A.7}
\end{eqnarray}%
where $2\leq i\leq n,$ $\nu =\frac{\xi }{\left\vert \xi \right\vert }.\ $In
addition, we have%
\begin{eqnarray}
&&\partial _{\xi ^{i}}\partial _{\xi ^{j}}\left( \left\vert \xi \right\vert 
\overset{.}{z}^{k}\left( x,\frac{\xi }{\left\vert \xi \right\vert },t\right)
\right) =\frac{g_{ij}}{\left\vert \xi \right\vert }\overset{.}{z}^{k}-\frac{1%
}{\left\vert \xi \right\vert ^{3}}\left( \underset{r=2}{\overset{n}{\sum }}%
g_{jr}\xi ^{r}\right) \left( \underset{r=2}{\overset{n}{\sum }}g_{ir}\xi
^{r}\right) \overset{.}{z}^{k}  \notag \\
&&-\frac{1}{\left\vert \xi \right\vert ^{2}}\left( \underset{r=2}{\overset{n}%
{\sum }}g_{jr}\xi ^{r}\right) \left( \overset{n}{\underset{s=1}{\sum }}%
\left( \partial _{\nu ^{s}}\overset{.}{z}^{k}\right) \frac{1}{\left\vert \xi
\right\vert ^{2}}\left( \xi ^{s}\underset{r=2}{\overset{n}{\sum }}g_{ir}\xi
^{r}\right) -\partial _{\nu ^{i}}\overset{.}{z}^{k}\right)  \notag \\
&&-\left( \partial _{\nu ^{i}}\overset{.}{z}^{k}\right) \frac{1}{\left\vert
\xi \right\vert ^{2}}\underset{r=2}{(\overset{n}{\sum }}g_{jr}\xi ^{r})-%
\overset{n}{\underset{s=1}{\sum }}\left( \partial _{\nu ^{s}}\overset{.}{z}%
^{k}\right) \frac{\xi ^{s}g_{ij}}{\left\vert \xi \right\vert ^{2}}  \notag \\
&&+\frac{2}{\left\vert \xi \right\vert ^{4}}\left( \underset{r=2}{\overset{n}%
{\sum }}g_{jr}\xi ^{r}\right) \left( \overset{n}{\underset{s=1}{\sum }}%
\left( \partial _{\nu ^{s}}\overset{.}{z}^{k}\right) \xi ^{s}\left( \underset%
{r=2}{\overset{n}{\sum }}g_{jr}\xi ^{r}\right) \right)  \notag \\
&&+\overset{n}{\underset{s=1}{\sum }}\frac{1}{\left\vert \xi \right\vert ^{3}%
}\xi ^{s}\underset{r=2}{(\overset{n}{\sum }}g_{ir}\xi ^{r})\left( \overset{n}%
{\underset{m=1}{\sum }}\left( \partial _{\nu ^{m}}\partial _{\nu ^{s}}%
\overset{.}{z}^{k}\right) \frac{1}{\left\vert \xi \right\vert ^{2}}\xi ^{m}%
\underset{r=2}{\overset{n}{\sum }}g_{ir}\xi ^{r}-\partial _{\nu
^{s}}\partial _{\nu ^{i}}\overset{.}{z}^{k}\right)  \notag \\
&&-\frac{1}{\left\vert \xi \right\vert ^{3}}\left( \underset{r=2}{\overset{n}%
{\sum }}g_{ir}\xi ^{r}\right) \overset{n}{\underset{m=1}{\sum }}\xi
^{m}\partial _{\nu ^{m}}\partial _{\nu ^{j}}\overset{.}{z}^{k}+\frac{1}{%
\left\vert \xi \right\vert }\partial _{\nu ^{j}}\partial _{\nu ^{i}}\overset{%
.}{z}^{k},  \TCItag{A.8}
\end{eqnarray}%
where $2\leq i,$ $j\leq n.$ Since $\overset{.}{z}^{k}\in C^{5}\left( \Omega
\left( d_{0}\right) \right) $ and the set $G^{\prime }$ is bounded and
closed, by (A.3)-(A.8), it is easy to see that%
\begin{equation}
\left\vert \left\vert \xi \right\vert ^{\left\vert \alpha \right\vert
-1}\partial _{\xi }^{\alpha }\left( \left\vert \xi \right\vert \overset{.}{z}%
^{k}\right) \right\vert \leq K_{2},\ 1\leq \left\vert \alpha \right\vert
\leq 4,  \tag{A.9}
\end{equation}%
in $\Omega ,$ where\ $K_{2}>0$ is independent of $(x,\xi )\in (D\times G)$.
In (A.9),$\ K_{2}$ depends on the norm of the vector function $\overset{.}{z}%
\left( x,\nu ,t\right) $ in the space $C^{5}\left( \Omega \left(
d_{0}\right) \right) $ and on the Euclidean distance between $G^{\prime }$
and$\ 0\in 
%TCIMACRO{\U{211d} }%
%BeginExpansion
\mathbb{R}
%EndExpansion
^{n-1}$ $(0\notin G^{\prime })$\ and the Euclidean diameter of $G^{\prime }$%
. Moreover, the following inequalities are valid in $\Omega $:%
\begin{equation}
\left\vert \left\vert \xi \right\vert ^{\left\vert \alpha \right\vert
+1}\partial _{\xi }^{\alpha }\left( \frac{1}{\left\vert \xi \right\vert }%
\right) \right\vert \leq K_{3},\ \left\vert \left\vert \xi \right\vert
^{\left\vert \alpha \right\vert }\partial _{\xi }^{\alpha }\left(
b(z)\right) \right\vert \leq K_{4},  \tag{A.10}
\end{equation}%
for$\ 0\leq \left\vert \alpha \right\vert \leq 4$, where\ $K_{3},$ $K_{4}$
depend on the same parameters as $K_{2}$ in (A.9).

Consequently, by the differentiation of an integral with respect to a
parameter, the boundedness of the functions $\left\vert \xi \right\vert 
\overset{.}{z}^{k},$ $\partial _{\xi ^{j}}\left( \left\vert \xi \right\vert 
\overset{.}{z}^{k}\right) $ and relations (A.9) - (A.10), we complete the
proof of the second assertion of Lemma 1.

\textbf{Proof of Lemma 2. }By (ii-c) of Lemma 1, we have%
\begin{equation*}
\xi ^{1}\partial _{\xi ^{1}}^{r}\left( \partial _{\xi ^{\prime }}^{\beta
^{\prime }}I\right) ,\ \xi ^{1}\partial _{\xi ^{1}}^{r}\left( \partial _{\xi
^{\prime }}^{\beta ^{\prime }}\partial _{x^{_{j}}}I\right) \in L_{1}\left( 
%TCIMACRO{\U{211d} }%
%BeginExpansion
\mathbb{R}
%EndExpansion
_{\xi ^{1}}^{1}\right) \cap L_{2}\left( 
%TCIMACRO{\U{211d} }%
%BeginExpansion
\mathbb{R}
%EndExpansion
_{\xi ^{1}}^{1}\right) ,\quad r+|\beta ^{\prime }|=4,\,0\leq r\leq 4
\end{equation*}%
for fixed $x\in D$ and $\xi ^{\prime }\in G^{\prime }$ $\left( \xi ^{\prime
}\neq 0\right) $. Then, in view of (i) of Lemma 1, (A.3)-(A.8) and the
properties of the Fourier transform, we have%
\begin{equation}
\eta ^{r}\partial _{\xi ^{\prime }}^{\beta ^{\prime }}\partial _{\eta }%
\widehat{I},\eta ^{r}\partial _{\xi ^{\prime }}^{\beta ^{\prime }}\partial
_{\eta }\partial _{x^{_{j}}}\widehat{I}\in C\left( D\times 
%TCIMACRO{\U{211d} }%
%BeginExpansion
\mathbb{R}
%EndExpansion
_{\eta }^{1}\times G^{\prime }\right) \cap L_{2}\left( 
%TCIMACRO{\U{211d} }%
%BeginExpansion
\mathbb{R}
%EndExpansion
_{\eta }^{1}\right) ,  \tag{A.11}
\end{equation}%
where$\ r+\left\vert \beta ^{\prime }\right\vert =4,\ 0\leq r\leq 4,$\ that
is, (iii) of Lemma 2 is proved.

Assertion (ii) of Lemma 2 is proved by using (i) and (ii-b) of Lemma 1.

By (i) and (ii-c) of Lemma 1, we have%
\begin{equation*}
\partial _{\xi ^{1}}^{r}\left( \partial _{\xi ^{\prime }}^{\beta ^{\prime
}}I\right) ,\ \partial _{\xi ^{1}}^{r}\left( \partial _{\xi ^{\prime
}}^{\beta ^{\prime }}\partial _{x^{_{j}}}I\right) \in L_{1}\left( 
%TCIMACRO{\U{211d} }%
%BeginExpansion
\mathbb{R}
%EndExpansion
_{\xi ^{1}}^{1}\right) \cap L_{2}\left( 
%TCIMACRO{\U{211d} }%
%BeginExpansion
\mathbb{R}
%EndExpansion
_{\xi ^{1}}^{1}\right)
\end{equation*}%
for $r+\left\vert \beta ^{\prime }\right\vert =4$ and $0\leq r\leq 4$, $x\in
D$ and $\xi ^{\prime }\neq 0,\in G^{\prime }$. The last relations and (i) of
Lemma 1 show that%
\begin{equation}
\eta ^{r}\partial _{\xi ^{\prime }}^{\beta ^{\prime }}\widehat{I},\ \eta
^{r}\partial _{\xi ^{\prime }}^{\beta ^{\prime }}\partial _{x^{_{j}}}%
\widehat{I}\in C\left( D\times 
%TCIMACRO{\U{211d} }%
%BeginExpansion
\mathbb{R}
%EndExpansion
_{\eta }^{1}\times G^{\prime }\right) \cap L_{2}\left( 
%TCIMACRO{\U{211d} }%
%BeginExpansion
\mathbb{R}
%EndExpansion
_{\eta }^{1}\right) ,  \tag{A.12}
\end{equation}%
consequently,%
\begin{equation}
\partial _{\xi ^{\prime }}^{\beta ^{\prime }}\widehat{I},\ \partial _{\xi
^{\prime }}^{\beta ^{\prime }}\partial _{x^{_{j}}}\widehat{I}\in C\left(
D\times \Delta _{\eta }^{\rho }\times G^{\prime }\right) ,  \tag{A.13}
\end{equation}%
for $\left\vert \beta ^{\prime }\right\vert \leq 4$. The first assertion of
Lemma 2 is ensured by (ii-a) of Lemma 1, (A.12) and (A.13).

\textbf{Proof of Lemma 3}. We examine equation (4.6) as a differential
equation for the function $\partial _{\eta }q$. Then, in view of the
characteristics, we can rewrite (4.6) as 
\begin{equation}
\frac{d}{ds}x^{1}=1,\ \frac{d}{ds}\xi ^{k}=-2\underset{j=2}{\overset{n}{%
\mathop{\textstyle \sum }}}\Gamma _{1j}^{k}\xi ^{j},\ \frac{d}{ds}\left(
\partial _{\eta }q\right) =\digamma _{2},\ k=2,3,...,n.  \tag{A.14}
\end{equation}%
By Lemma 2 we have%
\begin{equation}
\partial _{\xi ^{\prime }}^{\beta ^{\prime }}\digamma _{2},\text{ }\partial
_{\xi ^{\prime }}^{\beta ^{\prime }}\partial _{x_{j}}\digamma _{2}\in
C\left( D\times \Delta _{\eta }^{1}\times G^{\prime }\right) \cap
L_{2}\left( \Delta _{\eta }^{1}\right) ,  \tag{A.15}
\end{equation}%
for $0\leq \left\vert \beta ^{\prime }\right\vert \leq 2$, $1\leq j\leq n$.
Thus it follows from equalities (4.7), (A.14) and (A.15) that%
\begin{equation}
\partial _{\eta }q(x,\eta ,\xi ^{\prime })=\int_{x_{0}^{1}}^{x^{1}}\digamma
_{2}\left( \tau ,x^{\prime },\eta ,\zeta ^{\prime }(\tau ,\xi ^{\prime
})\right) d\tau ,  \tag{A.16}
\end{equation}%
where $x_{0}^{1}$ is the first component of the boundary point $%
(x_{0}^{1},x^{\prime })\in \partial D$ and the first component of $%
(x^{1},x^{\prime })\in D$ satisfies $x^{1}>x_{0}^{1}$. In (A.16), the
components of the vector $\zeta ^{\prime }(\tau ,\xi ^{\prime })=(\zeta
^{2}(\tau ,\xi ^{\prime }),...,\zeta ^{n}(\tau ,\xi ^{\prime }))$ satisfy
the system of differential equations%
\begin{equation*}
\frac{d}{d\tau }\zeta ^{k}=-2\underset{j=2}{\overset{n}{\mathop{\textstyle
\sum }}}\Gamma _{1j}^{k}\zeta ^{j},\quad 2\leq k\leq n
\end{equation*}%
with the initial condition $\ \zeta ^{\prime }(x^{1})=\xi ^{\prime }.$
Moreover, the uniqueness of the solution of this Cauchy problem with the
condition $\zeta ^{\prime }(x^{1})=\xi ^{\prime }\neq 0$ implies that $\zeta
^{\prime }(\tau ,\xi ^{\prime })\neq 0,$ $\tau \in \left[ x_{0}^{1},\text{ }%
x_{0}^{1}+d_{0}\right] ,$ $x_{0}^{1}\leq x^{1}\leq x_{0}^{1}+d_{0}$.

Since the function $\digamma _{2}\left( x,\eta ,\xi ^{\prime }\right) $ is
zero outside $D$ and the straight lines in $%
%TCIMACRO{\U{211d} }%
%BeginExpansion
\mathbb{R}
%EndExpansion
_{x}^{n}$ which are paralel to the coordinate axis $ox^{1}$ are geodesics of
the metric $\left( g_{ij}\right) ,$ the integral in (A.16) is taken on the
finite interval $(x_{0}^{1},x_{0}^{1}+d_{0}),$ where $d_{0}$ is the diameter
of the bounded domain $D$.\ On the other hand, by (4.2), (A.9), (A.10) and
Lemma 1, we see that the integral $\int_{-\infty }^{+\infty }u^{2}(x,\xi
)d\xi ^{1}$ converges uniformly with respect to the parameters $(x,\xi
^{\prime })\in D\times G^{\prime }$ and is continuous on $D\times G^{\prime
}.$\ Then by the Plancherel equality%
\begin{equation*}
2\pi \int_{-\infty }^{+\infty }u^{2}(x,\xi )d\xi ^{1}=\int_{-\infty
}^{+\infty }\left\vert \widehat{u}(x,\eta ,\xi ^{\prime })\right\vert
^{2}d\eta ,
\end{equation*}%
the integrals%
\begin{equation*}
\int_{0}^{+\infty }q^{2}(x,\eta ,\xi ^{\prime })d\eta ,\text{ }%
\int_{0}^{+\infty }p^{2}(x,\eta ,\xi ^{\prime })d\eta
\end{equation*}%
are continuous on $D\times G^{\prime }$.

In addition, by (4.2), (A.9) and (A.10), the integrals%
\begin{equation*}
\int_{-\infty }^{+\infty }(\partial _{\xi }^{\beta }u(x,\xi ))^{2}d\xi ^{1},%
\text{ }\int_{-\infty }^{+\infty }(\partial _{\xi }^{\beta }\partial
_{x^{_{j}}}u(x,\xi ))^{2}d\xi ^{1}
\end{equation*}%
converge uniformly with respect to the parameters $(x,\xi ^{\prime })\in
D\times G^{\prime }$ and are continuous on $D\times G^{\prime }$ for $%
\left\vert \beta \right\vert \leq 3,$ $1\leq j\leq n.$ Then, by similar
arguments as above one can prove that the integrals%
\begin{equation*}
\int_{0}^{+\infty }(\partial _{\xi ^{\prime }}^{\beta ^{\prime }}q(x,\eta
,\xi ^{\prime }))^{2}d\eta ,\text{ }\int_{0}^{+\infty }(\partial _{\xi
^{\prime }}^{\beta ^{\prime }}p(x,\eta ,\xi ^{\prime }))^{2}d\eta ,
\end{equation*}%
\begin{equation*}
\int_{0}^{+\infty }(\partial _{\xi ^{\prime }}^{\beta ^{\prime }}\partial
_{x^{j}}q(x,\eta ,\xi ^{\prime }))^{2}d\eta ,\text{ }\int_{0}^{+\infty
}(\partial _{\xi ^{\prime }}^{\beta ^{\prime }}\partial _{x^{j}}p(x,\eta
,\xi ^{\prime }))^{2}d\eta
\end{equation*}%
are continuous on $D\times G^{\prime }$ for $\left\vert \beta ^{\prime
}\right\vert \leq 3$. Therefore the integrals%
\begin{equation*}
\int_{0}^{+\infty }(\partial _{\xi ^{\prime }}^{\beta ^{\prime }}\digamma
_{2}(x,\eta ,\xi ^{\prime }))^{2}d\eta ,\ \int_{0}^{+\infty }(\partial _{\xi
^{\prime }}^{\beta ^{\prime }}\partial _{x^{_{j}}}\digamma _{2}(x,\eta ,\xi
^{\prime }))^{2}d\eta
\end{equation*}%
are continuous on $D\times G^{\prime }$ for $\left\vert \beta ^{\prime
}\right\vert \leq 2.$

On the other hand, by (A.16), we have%
\begin{equation}
q_{\eta }^{2}\leq (x^{1}-x_{0}^{1})\int_{x_{0}^{1}}^{x^{1}}\digamma
_{2}^{2}\left( \tau ,x^{\prime },\eta ,\zeta ^{\prime }(\tau ,\xi ^{\prime
})\right) d\tau .  \tag{A.17}
\end{equation}%
Since the function%
\begin{equation*}
\int_{0}^{+\infty }\digamma _{2}^{2}\left( x,\eta ,\xi ^{\prime }\right)
d\eta
\end{equation*}%
is continuous with respect to the parameters $(x,\xi ^{\prime })\in D\times
G^{\prime }$ and $D$, $G^{\prime }$ are closed and bounded, we can conclude
that there exists an integral%
\begin{equation*}
\int_{x_{0}^{1}}^{x^{1}}\left( \int_{0}^{+\infty }\digamma _{2}^{2}\left(
\tau ,x^{\prime },\eta ,\zeta ^{\prime }(\tau ,\xi ^{\prime })\right) d\eta
\right) d\tau
\end{equation*}%
which is bounded by a number $M>0$ which is independent of $(x,\xi ^{\prime
})$. Then, from (A.17) by the Fubini-Tonelli theorem, we have%
\begin{eqnarray*}
\int_{0}^{N}\left( \partial _{\eta }q(x,\eta ,\xi ^{\prime })\right)
^{2}d\eta &\leq &d_{0}\int_{0}^{N}\left( \int_{x_{0}^{1}}^{x^{1}}\digamma
_{2}^{2}\left( \tau ,x^{\prime },\eta ,\zeta ^{\prime }(\tau ,\xi ^{\prime
})\right) d\tau \right) d\eta \\
&=&d_{0}\int_{x_{0}^{1}}^{x^{1}}\left( \int_{0}^{N}\digamma _{2}^{2}\left(
\tau ,x^{\prime },\eta ,\zeta ^{\prime }(\tau ,\xi ^{\prime })\right) d\eta
\right) d\tau \\
&\leq &d_{0}M
\end{eqnarray*}%
for each $N>0$. By analogous reasoning, from (A.16), one can prove that%
\begin{equation}
\int_{0}^{N}(\partial _{\xi ^{\prime }}^{\beta ^{\prime }}\partial _{\eta
}q(x,\eta ,\xi ^{\prime }))^{2}d\eta ,\text{ }\int_{0}^{N}(\partial _{\xi
^{\prime }}^{\beta ^{\prime }}\partial _{x^{_{j}}}\partial _{\eta }q(x,\eta
,\xi ^{\prime }))^{2}d\eta \leq d_{0}M_{1}  \tag{A.18}
\end{equation}%
for $\left\vert \beta ^{\prime }\right\vert \leq 2$, where $M_{1}>0$ is the
maximum of the continuous functions%
\begin{equation*}
\int_{x_{0}^{1}}^{x^{1}}(\int_{0}^{+\infty }(\partial _{\xi ^{\prime
}}^{\beta ^{\prime }}\digamma _{2}(\tau ,x^{\prime },\eta ,\zeta ^{\prime
}(\tau ,\xi ^{\prime })))^{2}d\eta )d\tau
\end{equation*}%
and%
\begin{equation*}
\int_{x_{0}^{1}}^{x^{1}}(\int_{0}^{+\infty }(\partial _{\xi ^{\prime
}}^{\beta ^{\prime }}\partial _{x^{_{j}}}\digamma _{2}(\tau ,x^{\prime
},\eta ,\zeta ^{\prime }(\tau ,\xi ^{\prime })))^{2}d\eta )d\tau
\end{equation*}%
on $D\times G^{\prime }$. Inequalities (A.18) imply that $\partial _{\xi
^{\prime }}^{\beta ^{\prime }}\partial _{\eta }q,$ $\partial _{\xi ^{\prime
}}^{\beta ^{\prime }}\partial _{x^{_{j}}}\partial _{\eta }q\in L_{2}\left(
\Delta _{\eta }^{1}\right) $ for $\left\vert \beta ^{\prime }\right\vert
\leq 2$. Moreover, taking into account Remark 2 after Lemma 2, from (A.11)
we obtain%
\begin{equation*}
\partial _{\xi ^{\prime }}^{\beta ^{\prime }}\partial _{\eta }q(x,\eta ,\xi
^{\prime }),\text{ }\partial _{\xi ^{\prime }}^{\beta ^{\prime }}\partial
_{x^{_{j}}}\partial _{\eta }q(x,\eta ,\xi ^{\prime })\in L_{1}\left( \Delta
_{\eta }^{1}\right) \cap C\left( D\times \Delta _{\eta }^{1}\times G^{\prime
}\right) .
\end{equation*}%
In a similar manner, we can prove that%
\begin{equation*}
\partial _{\xi ^{\prime }}^{\beta ^{\prime }}\partial _{\eta }q(x,\eta ,\xi
^{\prime }),\text{ }\partial _{\xi ^{\prime }}^{\beta ^{\prime }}\partial
_{x^{_{j}}}\partial _{\eta }q(x,\eta ,\xi ^{\prime })\in L_{1}\left( \Delta
_{\eta }^{-1}\right) \cap C\left( D\times \Delta _{\eta }^{-1}\times
G^{\prime }\right) \text{.}
\end{equation*}%
Taking into account (4.5) and using arguments similar to the previous, it
can be shown that%
\begin{equation*}
\partial _{\xi ^{\prime }}^{\beta ^{\prime }}\partial _{\eta }p(x,\eta ,\xi
^{\prime }),\text{ }\partial _{\xi ^{\prime }}^{\beta ^{\prime }}\partial
_{x^{_{j}}}\partial _{\eta }p(x,\eta ,\xi ^{\prime })\in L_{1}\left( \Delta
_{\eta }^{\rho }\right) \cap L_{2}\left( \Delta _{\eta }^{\rho }\right) \cap
C\left( D\times \Delta _{\eta }^{\rho }\times G^{\prime }\right)
\end{equation*}%
for $\left\vert \beta ^{\prime }\right\vert \leq 2,$ $\rho =-1,1$.

Finally, from the foregoing relations we can write%
\begin{equation*}
\partial _{\xi ^{\prime }}^{\beta ^{\prime }}\partial _{\eta }\widehat{u},%
\text{ }\partial _{\xi ^{\prime }}^{\beta ^{\prime }}\partial
_{x^{j}}\partial _{\eta }\widehat{u}\in L_{1}\left( \Delta _{\eta }^{\rho
}\right) \cap L_{2}\left( \Delta _{\eta }^{\rho }\right) \cap C\left(
D\times \Delta _{\eta }^{\rho }\times G^{\prime }\right)
\end{equation*}%
for $\left\vert \beta ^{\prime }\right\vert \leq 2,$ $1\leq j\leq n.$ Then
the equality%
\begin{equation}
\partial _{\xi ^{\prime }}^{\beta ^{\prime }}q(x,\eta ,\xi ^{\prime
})=-\int_{\eta }^{\infty }\partial _{\xi ^{\prime }}^{\beta ^{\prime
}}\partial _{\tau }q(x,\tau ,\xi ^{\prime })d\tau  \tag{A.19}
\end{equation}%
holds, from which we have%
\begin{equation}
\partial _{\xi ^{\prime }}^{\beta ^{\prime }}q(x,+0,\xi ^{\prime
})=-\int_{0}^{\infty }\partial _{\xi ^{\prime }}^{\beta ^{\prime }}\partial
_{\tau }q(x,\tau ,\xi ^{\prime })d\tau  \tag{A.20}
\end{equation}%
for the point $(x,\xi ^{\prime })\in D\times G^{\prime }$. By virtue of
(A.18)-(A.20), we obtain%
\begin{eqnarray*}
\left\vert \partial _{\xi ^{\prime }}^{\beta ^{\prime }}q(x,\eta ,\xi
^{\prime })-\partial _{\xi ^{\prime }}^{\beta ^{\prime }}q(x,+0,\xi ^{\prime
})\right\vert ^{2} &=&\left\vert \int_{0}^{\eta }\partial _{\xi ^{\prime
}}^{\beta ^{\prime }}\partial _{\tau }q(x,\tau ,\xi ^{\prime })d\tau
\right\vert ^{2} \\
&\leq &\eta \int_{0}^{\eta }(\partial _{\xi ^{\prime }}^{\beta ^{\prime
}}\partial _{\tau }q(x,\tau ,\xi ^{\prime }))^{2}d\tau \\
&\leq &\eta d_{0}M_{1}.
\end{eqnarray*}%
Therefore, the function $\partial _{\xi ^{\prime }}^{\beta ^{\prime
}}q(x,\eta ,\xi ^{\prime })$ tends to $\partial _{\xi ^{\prime }}^{\beta
^{\prime }}q(x,+0,\xi ^{\prime })$ uniformly with respect to the parameters $%
(x,\xi ^{\prime })\in D\times G^{\prime }$ as$\ \eta \rightarrow +0$.
Consequently, $\partial _{\xi ^{\prime }}^{\beta ^{\prime }}q(x,+0,\xi
^{\prime })\in C\left( D\times G^{\prime }\right) $, since the function $%
\partial _{\xi ^{\prime }}^{\beta ^{\prime }}q(x,\eta ,\xi ^{\prime })\in
C\left( D\times G^{\prime }\right) $ for $\eta >0.$ Analogously, one can
prove that the functions $\partial _{\xi ^{\prime }}^{\beta ^{\prime
}}\partial _{x^{j}}q(x,+0,\xi ^{\prime }),$ $\partial _{\xi ^{\prime
}}^{\beta ^{\prime }}p(x,+0,\xi ^{\prime }),$ $\partial _{\xi ^{\prime
}}^{\beta ^{\prime }}\partial _{x^{j}}p(x,+0,\xi ^{\prime }),$ $\partial
_{\xi ^{\prime }}^{\beta ^{\prime }}q(x,-0,\xi ^{\prime }),$ $\partial _{\xi
^{\prime }}^{\beta ^{\prime }}\partial _{x^{j}}q(x,-0,\xi ^{\prime }),$ $%
\partial _{\xi ^{\prime }}^{\beta ^{\prime }}p(x,-0,\xi ^{\prime }),$ $%
\partial _{\xi ^{\prime }}^{\beta ^{\prime }}\partial _{x^{j}}p(x,-0,\xi
^{\prime })$ belong to the space $C\left( D\times G^{\prime }\right) $ for $%
\left\vert \beta ^{\prime }\right\vert \leq 2,$ $1\leq j\leq n$.

\bigskip

%The first author deceased in 2011. However, this work is mainly based on his
%ideas.

\textbf{Acknowledgements}. Most part of the paper has been written during
the stay of the second-named author at Department of Mathematical Sciences
of The University of Tokyo and the stay was supported by the program
"Leading Graduate Course for Frontiers of Mathematical Sciences and
Physics". The third-named author is supported by Grant-in-Aid for Scientific
Research (S) 15H05740 of Japan Society for the Promotion of Science.

\end{document}